\documentclass[review]{elsarticle}

\usepackage{lineno,hyperref,epstopdf}
\modulolinenumbers[5]

\usepackage{amssymb, amsmath, booktabs, color}
\newtheorem{lemma}{Lemma}
\newtheorem{theorem}{Theorem}%[section]
\newtheorem{remark}{Remark}
\biboptions{authoryear}
\begin{document}

\begin{frontmatter}

\title{A projection-based model checking for heterogeneous treatment effect}
\author[BNU]{Niwen Zhou}

\author[BNU]{{Xu Guo} }

\author[BNU,HKBU]{Lixing Zhu\corref{mycorrespondingauthor}\corref{thank} }
\cortext[thank]{
			The authors gratefully acknowledge two grants from the University Grants Council of Hong Kong  and a NSFC grant (NSFC11671042). \hspace{.2cm}}
\cortext[mycorrespondingauthor]{Corresponding author}
\ead{lzhu@hkbnu.edu.hk}

\address[BNU]{School of Statistics, Beijing Normal University, Beijing, China}

\address[HKBU]{Department of Mathematics, Hong Kong Baptist University, Hong Kong, China}

\begin{abstract}
In this paper, we investigate the hypothesis testing problem  that checks whether part of covariates / confounders significantly affect  the heterogeneous treatment effect given all covariates. This model checking is particularly useful in the case where there are many collected covariates  such that we can  possibly alleviate the typical curse of dimensionality.
 In the test construction procedure, we use a projection-based idea and a nonparametric estimation-based test procedure to construct an aggregated version over all projection directions. The resulting test statistic is then interestingly with no effect from slow convergence rate the nonparametric estimation usually suffers from. This feature makes the test behave like a global smoothing test to have ability to detect a broad class of local alternatives converging to the null at the  fastest possible rate in hypothesis testing. Also, the test can inherit the merit of lobal smoothing tests to be sensitive to oscillating alternative models.
The performance of the test is examined by  numerical studies and the analysis for a real data example for illustration.

\end{abstract}


\end{frontmatter}

%\linenumbers

\section{Introduction}

In this paper, we consider the testing problem for treatment effect model.
Let $D$ be the indicator variable of treatment and $Y$ the
outcome.  $D_i=0,1$ respectively means the $i$th individual
does not receive or receives treatment. The corresponding potential outcomes are then defined  as $Y_i(0)$ and $Y_i(1)$.  The observed
outcome can then be written as
$Y_i=D_iY_i(1)+(1-D_i)Y_i(0).$
An important quantity of interest in the literature is  the average treatment effect (ATE): $E(Y(1)-Y(0))$ %, the average treatment effects on the treated(ATT): $E(Y(1)-Y(0)\mid D=1)$
[see, e.g. \cite{Rosenbaum1983}, \cite{Hahn1998}]. To check the heterogeneity of ATE over a set $W$ of collected covariates, conditional (or heterogeneous) average treatment effect $E(Y(1)-Y(0)|W)$ (CATE) has been investigated.
 Note that CATE can capture the heterogeneity of a
treatment effect across the subpopulations defined by the index set of $W$. $CATE(W)$ is also called a contrast function in the precision medicine literature such as  \cite{shi2019}, which plays an important role in estimating optimal individualized treatment regime. In order to identify  this  function, other than the  common support assumption,  the unconfoundedness assumption is very important:
\begin{itemize}
 \item Assumption~1 (Unconfoundedness): $(Y(0),Y(1))\perp D\mid W$.
% \item Assumption~2(Common support):
% For some very small $c>0$,
% $c<p(W)<1-c$.
 \end{itemize}
Here $W=(X,Z)$ with $X$ and $Z$ being respectively $p-$ and $q-$dimensional vectors of covariates, and {$\perp$ stands for statistical independence.}
To make the paper self-contained, we write down the common support assumption as follows:
\begin{itemize}
 \item Assumption~2 (Common support):  For some very small $c>0$,
 $c<\pi(W)<1-c$
 \end{itemize}
{where $\pi(W)=E(D\mid W)$ is the propensity score function. }

 Based on these assumptions and others, most of existing estimation methods are for  $CATE(W)$ conditional on  all covariates $W$. See e.g. \cite{crump2008}, \cite{abrevaya2015}, \cite{hsu2017} and \cite{wager2018}. However, it may be the case that only the subset $X$ of $W$ is significantly useful for the average treatment effect such that $CATE(W)=CATE(X)$ in this case. This can then very much alleviate the curse of dimensionality in estimation and other further statistical analyses.   Such a dimension reduction structure needs an accompany  of model checking to prevent the working model possibly too parsimonious to lose some important covariates. More specifically, the hypotheses are:
 \begin{eqnarray}\label{null}
\begin{split}
&&H_{0}:\,\,\,P(E[Y(1)-Y(0)|W]=E[Y(1)-Y(0)|X])=1,\\
&&H_{1}:\,\,\,P(E[Y(1)-Y(0)|W]= E[Y(1)-Y(0)|X])<1.
\end{split}
\end{eqnarray}
There are no tests available for this issue in the literature although some are relevant.  \cite{crump2008} focused on testing whether $CATE(W)$ equals  zero  or a given constant. \cite{chang2015} and \cite{hsu2017} respectively proposed tests for the null hypothesis that  $CATE(W)$ or $CATE(X)$ is non-negative for all values of covariates. However, these tests cannot be used for the above testing problem.
To the best of our knowledge, the research described herewith is the first attempt to handle such a problem in the literature.  Further, we consider the situation that the dimensions of $X$ and $Z$, that is, $p$ and $q$ are fixed, but $p$ could be much smaller than  $q$. Then, under $H_0$, CATE can be estimated only conditional on a much lower dimensional covariates vector, $X$.

{
 The following are two special cases of $H_0$.
\begin{itemize}
\item Example~1 (Treatment effect heterogeneity)
\, Set $X=\emptyset$. In this very special case, the above null hypothesis becomes
{$$H_{02}:\,\,\,P(E[Y(1)-Y(0)|W]=E[Y(1)-Y(0)])=1.$$}
That is, under the null hypothesis, the treatment effect will not change with the value of $W$, thus the treatment effect does not have heterogeneity across the subpopulation defined by the value of $W$. The rejection of $H_{02}$ implies the necessity of  estimating conditional treatment effect. {Note that this test has been discussed by \cite{crump2008}.}
\item Example~2 (Significant conditional treatment effect)
\, Set $X=\emptyset$ and consider $E[Y(1)-Y(0)]=0$. we are still interested in testing:
{$$H_{03}:\,\,\,P(E[Y(1)-Y(0)|W]=0)=1.$$}
That is, when receiving a treatment
has no effect on outcomes for the overall population, we want to check whether the treatment is still significant for some subpopulations.
\end{itemize}
%Here we give a simple example to illustrate the introduced three null hypotheses. Assume $Z$ is longevity and $X$ is sex.
%Then $H_{0}$ can help us check whether the sex has effect for the treatment difference given the longevity.
%While $H_{02}$ aims to answer whether  the treatment has the same effect on different age groups. For $H_{03}$, we wish to known  whether the program is significant for people at some specific age groups, even, on average, it does not affect outcomes.}
Since the test construction for $H_{02}$ and $H_{03}$ can be relatively easier in our methodology than that for  $H_{0}$, we then only deal with $H_0$ in the following. % and then give a brief discussion on the test constructions for the  other two testing problems.  %Clearly, $H_{0}$ is closely related with the typical significance testing problem. However, there is a big difference. For the latter testing problem, the response is generally observable. While for $H_{0}$, only $Y(0)$ or $Y(1)$ is observed for every individual. This puts a big challenge for the testing problem. Note that the typical significance testing problem with known response can be regarded as a special case of our work. Thus it can be conducted in the same paradigm, which will be described later.

{Clearly, we first need to  estimate the conditional mean functions $E(Y(1)-Y(0)\mid X)$ and $E(Y(1)-Y(0)\mid W)$.  As we do not assume any parametric model structure for these functions,  nonparametric estimation is applied.
As commented above,  when the dimension $p$ of the covariates $X$ is high, any  nonparametric estimation would be inefficient and thus has negative effect for the performance of constructed test.  We then review some typical  methods for regressions first. There are a number of proposals available in the literature, but we only name a few to comment on their pros and cons. For local smoothing tests for regressions in the literature, one of methods was proposed by  \cite{zheng1996}. %\cite{zheng1996consistent}
But it can only detect local alternatives distinct from the null at the rate of order $n^{-1/2}h^{-(p+q)/4}$ where $h$ is the bandwidth going to zero at a certain rate in nonparametric estimation. { This drawback can be found in other typical local smoothing test literature, see e.g. \cite{fan1996consistent}, \cite{zhang2004} and \cite{guo2016} although the rate could be proved, in some special model structures, to $n^{-1/2}h^{-1/4}$ when some dimension reduction approaches are applied. }Thus,  we also wish that it have the nice properties of global smoothing-based tests in the literature for regressions to detect local alternatives distinct from the null at the fastest possible rate of order $n^{-1/2}$ in hypothesis testing. {See \cite{stute1998}, \cite{zhu2003} and \cite{khmaladze2009} %\cite{stute1998bootstrap}
for such types of tests.}   %Another way is to use a direction $\beta$ that is randomly selected from all directions such that  test is still omnibus. %Further,  the test statistics based on typical kernel-based estimators can detect the local alternatives distinct from the null at a rate depended on the bandwidth $h$ and the dimension of $X$ or $W$, which is much slower than $n^{-1/2}$. These are also common drawbacks in typical nonparametric significance tests. See e.g.
%\cite{lavergne2000},  \cite{lavergne2015}.
%Therefore, these test statistics are asymptotically less powerful for detecting alternative models.
}To make a test sensitive, to a certain extent, to oscillating/high-frequency alternative models, it would be good to construct a test that is based on a local smoothing test structure by using some nonparametric estimation for the involved functions. On the other hand, to have the test more powerful to detect smooth local alternatives, we also wish it to have the features global smoothing tests share. Therefore, we combine two ideas to achieve these goals. First, to alleviate this dimensionality difficulty, we suggest a projection-based test that uses projected covariates $ \beta^{\top}W$. It is clear that we cannot simply use only one or a few projections to construct a test otherwise, it will be a directional test.  To make test omnibus against all alternatives, we then use the projected covariates at all projection directions in an  aggregation manner. From \cite{zhu1998}, \cite{escanciano2006}, \cite{stute2008}, and \cite{lavergne2012}, we anticipate that the dimensionality issue could be largely alleviated as for regressions, these tests can reach the rate much faster than $n^{-1/2}h^{-(p+q)/4}$.
Thus, all these tests can very much improve the performance  in high-dimensional scenarios. But these tests are either still typical nonparametric estimation-based local smoothing tests that can detect local alternatives at slower rate than $1/\sqrt n$ or typical empirical process-based global smoothing tests that are less sensitive to high-frequency alternative models. Taking this issue into consideration, we consider constructing a test that is based on local smoothing  technique and then is transferred to a final pairwise distance-based test. Under certain regularity conditions, the limiting null distribution of this test statistic can then be free of the nonparametric estimation for the conditional moment on the whole $W$ such that the test behaves like a global smoothing test and at the same time, shares the sensitivity to high-frequency models to certain extent. This will be demonstrated in the numerical studies. Another feature of the test is worthwhile to mention: although the function $E[Y(1)-Y(0)|X]$ under the null hypothesis indispensably requires nonparametric estimation, it does not make a slow-down of the resulting rate of convergence and the test can still share all features global smoothing tests have.

%Three features of the test are worthwhile to mention.
%{
%\begin{itemize}
%\item First, we introduce a new kind testing problem for treatment effect model, which is more general and more informative than previously ones in literature.
%\item Second, we propose an efficient testing procedure for this new testing problem, which is omnibus to detect all global alternative hypothesis and can achieve the fastest rate to  detect local alternatives hypothesis.
%\item Lastly, the pairwise distance part in the test statistic play an important role in inherits some features of local smoothing tests and make it possible to alleviate the curse of dimensionality.
%\end{itemize}
%}

{The rest of this paper is organized as follows. In Section 2, we describe the test statistic construction. The asymptotic properties of the test statistic under the null, global and local alternative hypothesis are investigated in Section 3. %We also display the asymptotic properties in the special cases in the two examples we discussed before.
In Section 4, we examine the finite sample performance of our test through simulations and  apply it to a real
data example for illustration in Section~5. Some conclusions are presented in Section 5, and the proofs of the
theoretical results are postponed to Appendix.}

\section{The  test statistic construction}

%Without loss generality, assume the random sample $(Y_i,W_i,D_i), 1\leq i\leq n$, from $(Y,W,D)$
%
%The propensity score, $E(D|X,Z)$ denotes as $p(X,Z)$.
%
%For identification and inference, we make the following two common assumptions:
%\begin{itemize}
% \item[(i)] Unconfoundedness: $(Y(0),Y(1))\perp D\mid X,Z$.
%
% \item[(ii)] Common support: For some very small $c>0$,
% $c<p(X,Z)<1-c$.
% \end{itemize}
%
%Let $W=(X,Z)$.

Note that under the unconfoundedness  assumption and common support assumption, the conditional treatment effect can be identified as:
\begin{eqnarray*} E[Y(1)-Y(0)|W]=E\bigg[\frac{DY}{\pi(W)}-\frac{(1-D)Y}{1-\pi(W)}\Big |W\bigg],
\,\,\,\,E[Y(1)-Y(0)|X]=E\bigg[\frac{DY}{\pi(W)}-\frac{(1-D)Y}{1-\pi(W)}\Big |X\bigg].\end{eqnarray*}
Let $  Y^*=\frac{DY}{\pi(W)}-\frac{(1-D)Y}{1-\pi(W)}$.
Then $H_{0}$ can be rewritten as follows:
\begin{eqnarray} H_{01}:\,\,\,P(E[Y^*|W]=E[Y^*|X])=1.\end{eqnarray}

% {\color{red}Then we extend to investigate the situation that $p$ is fixed, while $q$ is divergent.}

%{\color{blue} we want to asses the significance of $X\in R^q$ in the nonparametric regression of $Y\in R$ on $W\in R^p$ and $X$.}

%The observed data from n subjects can be summarized as
%$O_i = Si(1),A( i1),Si(2),A( i2),Yi, i = 1,...,n,$
%which are assumed to be independently and identically distributed copies of O. Let $Y^*\in R,~X\in R^{p},~Z\in R^{q}$ are random vectors.

Note that we consider the case where the propensity score is a function of $W$, rather than $X$. This is because the propensity score is a probability for treatment $D$ when the covariates are given. Thus, the decision on whether giving treatment is based on all covariates / confounders. While the testing problem is for treatment effect after giving the decision on treatment. Thus, this is a reasonable scenario.

Define $g(X)=E(Y^*\mid X)$ and $e=Y^*-g(X),~u=Y^*-E(Y^*\mid W).$ Thus, under the null hypothesis, $e=u$ with $E(e\mid W)=0$, otherwise $E(e\mid W) \neq 0$. Hence, it is reasonable to directly construct a test statistic based on the sample analogue of $E[eE(e\mid W)]=E\{[E(e\mid W)]^2\}\geq 0$  with the equality holds if and only if $E(e\mid W)=0$. This idea is similar to that in \cite{zheng1996}. However, in order to get the sample analogue of $E[eE(e\mid W)]$ without a model misspecification risk, a nonparametric estimation of $E(e\mid W)$ is required. Note that $W=(X^{\top},Z^{\top})^{\top}$ with $X\in R^p$ and $Z \in R^q$. Hence, any nonparametric estimation of $E(e\mid W)$  suffers from the curse of dimensionality when the dimension $q$ of possible insignificant variables $Z$ is large, even moderate.
 This motivates us to construct test statistic based on a method with projection directions. To this end, we first give a lemma about the equivalence between function with original covariates and that with projected covariates below.

% comes from the fact that employing a high dimensional nonparametric estimation for $E(e\mid W)$ will typically result in poor estimation accuracy and the curse of high dimensionality.

\begin{lemma}\label{lemma1}
  $E(Y^*\mid W)=E(Y^*\mid X)$ holds if and only if
$E(e\mid \alpha^{\top}W)=0$ holds for all $\alpha \in R^{p+q}.$ {Further the equality $E(e\mid \alpha^{\top}W)=0$ holds if and only if $\int \{[E(e\mid \alpha^{\top}W)]^2f_\alpha(\alpha^{\top}W)\}\mu(\alpha)d\alpha=0$ when the function $E(e\mid \alpha^{\top}W)$ is continuous about $\alpha$.}
\end{lemma}
{
 %This lemma can be deduced from \cite{bierens1982}.
 Similar conclusion can be also found in \cite{zhu1998}, \cite{escanciano2006}, \cite{lavergne2012} and \cite{li2019}.

Note that Lemma~1 implies that, under the null hypothesis,
\begin{eqnarray}
\int E\{[E(e\mid \alpha^{\top}W )]^2f_\alpha( \alpha^{\top}W )\}\mu(\alpha)d\alpha=0.
\end{eqnarray}
 While under the alternative hypothesis, there exist some $\alpha_*\in
R^{p}$ such that $E(e\mid {\alpha_*}^{\top}W
)\neq 0$ and by the continuity of this function with respect to  $\alpha$, there is a neighborhood $\alpha$ whose measure   is positive and $E(e\mid {\alpha}^{\top}W
)\neq 0$ for all $\alpha$ in the neighborhood. Thus, it follows that \begin{eqnarray}\label{basic_stat}
\int E\{[E(e\mid \alpha^{\top}W )]^2f_\alpha( \alpha^{\top}W )\}\mu(\alpha)d\alpha >0.
\end{eqnarray} Also note that $E\{[E(e\mid \alpha^{\top}W )]^2f_\alpha( \alpha^{\top}W )\}=E\{eE(e\mid  \alpha^{\top}W ) f_\alpha( \alpha^{\top}W )\}$.
The above argument implies that we can use the sample analogue of $\int  E\{eE(e\mid  \alpha^{\top}W ) f_\alpha( \alpha^{\top}W )\}\mu(\alpha)d\alpha$ to construct a test statistic.
}

When an independent and identically distributed (i.i.d.) random sample $\{(Y_i,D_i,X_i,Z_i)\}_{i=1}^n$ is available with a parametric  propensity score function $\pi(W, r_0)$ where $r_0$ is an unknown parameter vector of dimension $d$,  we first estimate $Y_i^*$ by $\hat Y_i^*=\left(\frac{D_i}{ \pi(W_i,\hat r)}-\frac{1-D_i}{1-  \pi(W_i,\hat r)}\right)Y_i$, where $\hat r$ in $\pi(W,\hat r)$ is a maximum likelihood estimator of $r_0$.
The test statistic is defined as
\begin{eqnarray}
\tilde T_n=\int\frac{1}{n(n-1)}\sum_{i=1}^n\sum_{j\neq i}\hat e_i\hat e_j  H_{h_1}
\left(\alpha^{\top}W_i-\alpha^{\top}W_j\right)\mu(\alpha)d\alpha.
\end{eqnarray}
This is the sample version of (\ref{basic_stat}) where  $
\hat E( e_j\mid \alpha^{\top}W_{i})=  \frac{1}{(n-1)  } \frac{\sum_{j\neq i}^n \hat e_j H_{h_1}
\left(\alpha^{\top}W_j-\alpha^{\top}W_i\right)}{\hat f_\alpha(\alpha^{\top}W_{i})}$ is the kernel estimator of $E(e\mid \alpha^{\top}W )$ {with $H(\cdot)$ being a kernel function, $H_{h_1}(\cdot)=H(\cdot/h_1)/h_1$ and $h_1$ is the bandwidth.
 }$\hat e_i=\hat Y_i^*-\hat g(X_i)$ with $\hat g(X_i)= \sum_{j\neq i}^n   w_{ij}\hat Y_j^*$, $w_{ij}=\frac{1}{(n-1) }  \mathcal{K}_h\left(X_j-X_i \right)\big/\hat f(X_i)$. The density estimator is  $\hat f(X_i)=\frac{1}{(n-1)} \sum_{j\neq i}^n  \mathcal{K}_h\left(X_j-X_i \right)$, {where $\mathcal{K}$ is a multivariate kernel function and $\mathcal{K}_h(\cdot)=\frac{1}{h^p}\mathcal{K}
\left(\frac{\cdot}{h}\right)$ with a bandwidth $h$.}

Although  nonparametric kernel estimation for $E(Y^*\mid X)$ is inevitable, this test statistic only involves the integral of univariate $\alpha^{\top}W$ over all $\alpha$ rather than the original high-dimensional $W$. Thus,  $\tilde T_n$ could greatly mitigate the dimensionality problem due to the nonparametric estimator of $E(e\mid \alpha^{\top}W)$. However, it can be expected that the asymptotic properties of $\tilde T_n$ will still be related to the bandwidth $h_1$ in a nonparametric estimation nature so that the convergence rate would be slower than $1/\sqrt{n}$.
To tackle this problem, we adapt the idea in \cite{li2019} to transform the nonparametric estimation based test  into the pairwise distance-based one so that the convergence rate can be free of the bandwidth parameter $h_1$. To be specific, let $H(u)=\frac{1}{\sqrt{2\pi}}\exp\left(-\frac{u^2}{2}\right)$ and consider $\alpha\sim N(0,h_1^2I_p),$ using the direct consequence of $Lemma~2$ in \cite{li2019}.  It can be shown that $\tilde T_n$ is proportional to $T_n$, i.e. $\tilde T_n=\frac{1}{h_1}T_n$ with
\begin{eqnarray}\label{tn}
\tilde T_n=\frac{1}{h_1n(n-1)}\sum_{i=1}^n\sum_{j\neq i}\hat e_i\hat e_j \mathcal{B}_{ij},
\end{eqnarray}
where $\mathcal{B}_{ij}=\frac{1}{\sqrt{1+d_{ij}}}$ with $d_{ij}=\|W_i-W_j\|^2.$ Here $\|\cdot\|$ denotes the Frobenius norm throughout this paper.
Note that $h_1$ is just a constant outside the sum in $\tilde T_n$, therefore, we can use  $T_n=h_1\tilde T_n$ that is then free of the bandwidth $h_1$.

%
%\begin{remark}
%{
%We can also get the same test statistic $T_n$ in (\ref{tn}) via a fourier transform characterization of conditional moment. This method is also prevalent in literature to mitigate the impact of dimensionality, see \cite{bierens1982} and \cite{guo2019}. To be specific, note that $E(e\mid W)=0$ is equivalent as $E(e \exp^{it^{\top }W})=0$ due to the uniqueness of a function's Fourier transform. Hence we can construct a test statistic based on the sample analogy of
%$T_\omega=\int |E(e \exp^{it^{\top }W})|^2\omega(t)dt$.
%Further, if we assume $t$ have a $(p+q)$ dimensional symmetric laplace density function mentioned in \cite{kozubowski2013}, i.e. $t \sim GAL_r(2^{1/(p+q)}\mathbb{I}_{p+q},0,\frac{1}{2})$, it can be deduced that \begin{eqnarray}
%T_\omega=  E\left(e e'\frac{1}{\sqrt{1+\|W-W'\|^2}}\right),
%\end{eqnarray} where $\mathbb{I}_{p+q}$ denotes the $(p+q)$ dimensional identity matrix. Thus if we construct the test statistic based on the sample analogy of $T_\omega$, we can also get $T_n$ in (\ref{tn}).
%
%Thus, no matter we construct our test statistic via the spirit of bridging idea or fourier transform characterization, $T_n$ can be a omnibus test statistic to detect whether the two conditional treatment effects, $E(Y(1)-Y(0)\mid X)$ and $E(Y(1)-Y(0)\mid W)$, is equivalent.
%
%}
%\end{remark}
\begin{remark}
{
It is worth mentioning that the constructed test could still inherit some features of existing local smoothing tests. Recall that  $\mathcal{B}_{ij}=\frac{1}{\sqrt{1+d_{ij}}}$ with $d_{ij}=\|W_i-W_j\|^2$. Thus, $T_n$ captures more information from closely related observations.
 This property ensures that no matter the alternatives are either highly frequent or lowly frequent, the test  $T_n$ could be  workable to detect them. This merit can be confirmed by the numeric studies below.
 }
\end{remark}
%\begin{remark}
%{
%We cannot deny that $T_n$ can mitigate the curse of dimensionality. However, when the dimension of $W$, $(p+q)$, tends to infinity, the parameter $d_{ij}=\sum_{k=1}^{p+q}(W_{ik}-W_{jk})^2$ cannot be ignorable, so that $T_n$ will tend to zero no matter under the null or alternative hypothesis. In other words, $T_n$ will be powerless when the dimension of $W$ diverge and a new test procedure need to be explored in future work to tackle this problem.    }
%\end{remark}

\section{Asymptotic properties}

In order to get the asymptotic behaviour of $T_n$, the following assumptions are required:
\begin{itemize}
\item Assumption~3(Sampling): The observation data, $\{(Y_i,D_i,X_i,Z_i)\}_{i=1}^n$, is an independent and identically distributed random sample of size $n$ from the joint distribution of the vector $(Y,D,X,Z)$.
\item Assumption~4(Distribution): the density of $X$, $f(x)$, is bounded away from zero and infinity,  $s$-times continuously differentiable on its support $\Omega$. $E(Y^*\mid X)=g(X)$ is continuously differentiable.
\item Assumption~5(Moments): $E(u^4)<\infty$, $E(\eta^4)<\infty$  and $E\|W\|^2<\infty$.
\item Assumption~6(Kernel): For $p$-dimensional $u$, $\mathcal{K}(u)$ is a bounded kernel that is symmetric around zero, and $s$ times continuously differentiable and  of order $s$: $\int \mathcal{K}(u)du=1$, $\int u_1^{p_1}\cdots u_p^{p_p}\mathcal{K}(u)du=0$ for all nonnegative integers $p_1,\cdots,p_p$ such that $1\leq \sum_{i}p_i<s, $ and  nonzero when $ \sum_{i}p_i=s$.

\item  Assumption~7(Bandwidth): $h \rightarrow 0$, $nh^{2p} \rightarrow \infty$ and $nh^{2s}\rightarrow 0$ as $n \rightarrow \infty$.
\item Assumption~8(Propensity score estimator):
 The propensity score $\pi(W)=E(D
 \mid W)$ has a parametric form $\pi(W,r_0)$, and the function $\pi(W,r)$ is bounded away from zero and  has bounded continuous partial derivatives up to order 2 with respect to $r\in \Theta \subset R^{d}$, $d<\infty$.
%The estimator $\hat r$ of the propensity score model $\pi(W,r)$, $r \in \Theta \subset R^{p+q}$, $p+q<\infty$, satisfies $\sup_{W}|\pi(W,\hat r)-\pi(W,r_0)|=O_p(n^{-1/2})$ where $r_0\in \Theta$ such that $E(D\mid W)=\pi(W,r_0)$.
\end{itemize}
Here $\eta=\bar Y-E(\bar Y\mid X)$ with $\bar Y=\left(\frac{-D}{\pi(W,r_0)^2}-\frac{1-D}
 {(1-\pi(W,r_0))^2}\right)Y\nabla \pi(W,r_0)$, { and $\nabla\pi(W,r_0)$ stands for the  partial derivative with respect to $r_0$.}

Assumptions $3\sim 5$ are commonly used to guarantee the asymptotic normality of the test statistic.
As our test statistic uses nonparametric estimation for $E(Y^*\mid X)$, assumptions $6\sim 7$ are designed to
ensure the nonparametric estimator well-behaved, which are also widely used in the nonparametric estimation literature.
The condition  on $\mathcal{K}$ is for convenience of theoretical analysis. The following proof can be extended to kernels with exponential tails. Assumption~8 is standard in the literature to obtain $\sqrt{n}$-consistent estimation for $r_0$ in $\pi(W,r_0)$, see e.g. \cite{yao2010} and \cite{lin2018}. Based on this assumption, we can get the following lemma.

{
\begin{lemma}
 Under Assumption~8, the maximum likelihood estimator $\hat r$ has the following asymptotically linear representation:
\begin{eqnarray}
\hat r-r_0=\frac{1}{n}\sum_{i=1}^n  R_i+o_p(n^{-1/2}).
\end{eqnarray}
and its asymptotic distribution is
\begin{eqnarray}
\sqrt{n}(\hat r-r_0)\stackrel{\mathcal{D}}{\longrightarrow}N(0,\Sigma^{-1}).
\end{eqnarray}
where $R_i=\Sigma^{-1}\frac{
\nabla\pi(W_i,r_0)(D_i-\pi(W_i,r_0))}{\pi(W_i,r_0)(1-\pi(W_i,r_0))}$, $\Sigma=E\left[\frac{\nabla\pi(W,r_0)\nabla\pi(W,r_0)^{\top}}
{\pi(W_i,r_0)(1-\pi(W_i,r_0))}\right].$
\end{lemma}
}
This lemma can be found in \cite{yao2010}. Based on this lemma, we return to investigating the asymptotic distribution of $T_n$ under the null and alternative hypothsis.
%Kernels satisfying Assumption 7 can be constructed by taking products of higher order univariate kernels. A general method for obtaining higher order kernels from ¡°regular¡± ones is described, for example, by Imbens and Ridder (2009).
%expositional convenience only; we can extend the proof of Theorem 1 to kernels with exponential tails

\subsection{Asymptotic behavior under the null hypothesis}

%%%%%%%%%%%
% Introduce some notations to start the theory dis
%%%%%%%%%%%
%{\color{red}
%Before introducing asymptotic theories of $T_n$, we first define some important related notations. Let $H_j=H(W_j)$.
%}

We have the following asymptotic results under the null hypothesis.
\begin{theorem}\label{null_them}
Under Assumptions $1\sim 8$ and  $E[Y(1)-Y(0)\mid W]=E[Y(1)-Y(0)\mid X]$ holds with probability 1, the test statistic $T_n$ in (\ref{tn}) satisfies
\begin{eqnarray}
nT_n\stackrel{\mathcal{D}}{\longrightarrow}
\sum_{i=1}^\infty\lambda_i(Z_i^2-1)+2\nu_1^{\top}\nu_2+
\nu_1^{\top}A\nu_1+\mu^*
\end{eqnarray}
Here {$\nu_1\sim N(0,\Sigma^{-1})$, $\nu_2\sim N(0,\Sigma_1)$  with $\Sigma_1=E(e_1^2 \bar H_1\bar H_1^{\top})$ and $\bar H_1=E(\eta_2{\mathcal{B}}_{12}\mid W_1)-\frac{1}{2}E\left[(\eta_2\bar w_{31}
+\eta_3\bar w_{21})\mathcal{B}_{32}\mid W_1\right]$.}
$\mu^*=-E[(e_1^2\bar w_{12}+e_2^2\bar w_{21})
\mathcal{B}_{12}]+\frac{1}{3}E(e_1^2\bar w_{21}
\bar w_{31}\mathcal{B}_{23}
+e_2^2\bar w_{12}\bar w_{32}\mathcal{B}_{13}
+e_3^2\bar w_{12}\bar w_{23}\mathcal{B}_{12})$ with $\bar w_{ij}=  \mathcal{K}_h\left(X_j-X_i \right)\big/f(X_i)$. And $A=E(\eta_1\eta_2^{\top}\mathcal{B}_{12})$.

Let $Z_i$'s be independent standard normal random variables and $\lambda_{i}$'s  the eigenvalues of the integral equation
\begin{eqnarray}
\int L(\chi_1,\chi_2)\phi_i(\chi_2)dF(\chi_2)=\lambda_i\phi_i(\chi_1)
\end{eqnarray}
 with $\chi_i=(e_i,W_i)$ and $\phi_i(\chi)$ being the associated orthonomal eigenfunctions.
 Let $$\tilde{\mathcal{B}}_{ij}=\frac{1}{2}
 E\left(\bar w_{ti}
\mathcal{B}_{jt}+\bar w_{tj}\mathcal{B}_{it}\mid W_i,W_j\right)$$
 and $$\bar{\mathcal{B}}_{ij}=\frac{1}{2}
 E\left(\bar w_{ti}\bar w_{kj}\mathcal{B}_{tk}
 +\bar w_{ki}\bar w_{tj}\mathcal{B}_{kt}\mid W_i,W_j\right).$$
 Write $L(\chi_i,\chi_j)=e_ie_j(\mathcal{B}_{ij}
 -2\tilde{\mathcal{B}}_{ij}+\bar{\mathcal{B}}_{ij})$ \end{theorem}

 {
 Obviously, under the null hypothesis $H_0$, $T_n=O_p\left(\frac{1}{n}\right)$ implies the fact that $T_n$ converges to zero very quickly when $H_0$ is true. This will lead to a sensitive test to detect local alternatives close to the null at a fastest possible rate in hypothesis testing.  As $T_n$ contains the nonparametric estimation of $g(X)$ and the parametric propensity score estimation, the limiting null distribution of $T_n$ is  intractable. Thus we use the wild bootstrap  to approximate the null distribution. See  the details in Section \ref{num_study}.

}

\subsection{Power study}
{
 To examine the power performance of $T_n$, we consider the following sequence of local alternative hypotheses as:
\begin{eqnarray}\label{h1n}
H_{1n}:~P(E(Y(1)-Y(0)\mid W)=g(X)+a_nH(W))=1,~W=(X,Z).
\end{eqnarray}
Recall that $g(X)=E(Y(1)-Y(0)\mid X)$. Thus fixed $a_n$ corresponds to the global alternative model and when $a_n$ goes to zero, the sequence is about the local alternative hypotheses.
To smooth the theoretical analysis, the following assumption is added:
\begin{itemize}
\item Assumption~9(Alternatives):  $E(H(W)^2)<\infty$ and $E[H(W)\mid X]=0$.
\end{itemize}
The moment condition of $H(W)$ is  commonly assumed and the condition $E[H(W)\mid X]=0$ can be also found in \cite{lavergne2015}. {This condition is parallel to the unconditional one in \cite{lavergne2000} that is $H(W)\equiv 0$ when $a_n=0$.}}  we have the following theorem.

\begin{theorem}\label{alter_them}
Suppose Assumptions $1\sim 9$ hold. Then under the local alternative hypotheses in (\ref{h1n}), the following results can be obtained:
\begin{itemize}
\item[(1)] Under the global alternative hypothesis $H_{1n}$ with a fixed $a_n>0$,
\begin{eqnarray}
T_n\stackrel{P}{\longrightarrow} \mu>0.
\end{eqnarray}
Here $\mu=a_n^2E(H_1H_2\mathcal{B}_{12}).$
\item[(2)] Under the local alternative hypothesis $H_{1n}$ with $\sqrt{n}a_n\rightarrow \infty$,
\begin{eqnarray}
T_n/a_n^2\stackrel{P}{\longrightarrow} \mu_0>0.
\end{eqnarray}
Here $\mu_0=E(H_1H_2\mathcal{B}_{12}).$
\item[(3)] Under the local alternative hypothesis $H_{1n}$ with $a_n=n^{-1/2}$,
\begin{eqnarray}
nT_n\stackrel{\mathcal{D}}{\longrightarrow}\sum_{i=1}^
\infty\lambda_i(Z_i^2-1)+2\nu_1^{\top}\nu_2+
\nu_1^{\top}A\nu_1+N(\tilde \mu,\sigma^2).
\end{eqnarray}
\end{itemize}
$\tilde \mu=\mu_0+\mu^*,$ $\sigma^2=E(u_1^2\tilde H_1^2)+\alpha^{\top}\Sigma^{-1}\alpha$. Here $\tilde H_1=2E(H_2\mathcal{B}_{12}\mid W_1)+
\left[(H_2\bar w_{31}+H_3\bar w_{23})
\mathcal{B}_{23}\mid W_1\right]$ and $\alpha=E(H_1\eta_2\mathcal{B}_{12}).$
\end{theorem}

This theorem implies that when the local alternatives converge to the null hypothesis at a slower rate $a_n=O(n^{-c})$ than $O(n^{-1/2})$ for  $0\le c<1/2$, $nT_n\rightarrow \infty$ in probability at the rate of $n^{1/2-c}$. Thus the test is consistent as its asymptotic power tends to 1. Further, the test $T_n$ can still detect the local alternatives that are  distinct from the null hypothesis at a fastest possible rate $\sqrt{n}$. This is the typical feature existing global smoothing tests in the literature for regressions share.
\begin{remark}
From Theorem~{\ref{null_them}}, we can see that its limiting null distribution is rather complicated in formula. It is partly because of the effect from the estimation of  propensity score function. But it is interesting that the nonparametric estimation does not cause a slowdown of the convergence rate of $T_n$ to its weak limit.

%What's more, our framework can be also suitable for
%$$H_{04}:\,\,\,E[D|X,Z]=E(D|X).$$

\end{remark}

\section{Numerical studies}\label{num_study}
In this section we carry out two sets of simulation studies to evaluate the performance of the proposed test under different model settings. As the limiting null distribution of the proposed test is intractable,  we use the wild bootstrap approximation to determine  critical values. The procedure is given as follows:
\begin{itemize}
\item[(1)] For a given random sample $\{ (Y_i,W_i,D_i):i=1,\cdots,n\}$,  obtain $\hat e_i=\hat Y_i^*-\hat g(X_i)$ with $\hat Y_i^*=\left(\frac{D_i}{ \pi(W_i,\hat r)}-\frac{1-D_i}{1-  \pi(W_i,\hat r)}\right)Y_i$ and $\hat g(X_i)=\frac{1}{(n-1)} \sum_{j\neq i}^n  w_{ij}\hat Y_j^*$. Here $\hat r$ is the maximum likelihood estimator of $r_0$ in $\pi(W,r_0)$.
\item[(2)] Generate a bootstrap sample $\{ (\tilde Y_i^*,W_i,D_i):i=1,\cdots,n\}$ with the outcome variables as $\tilde Y_i^*=\hat g(X_i)+\tilde e_i$. Here $\tilde e_i=\zeta_i\hat e_i$ and $\zeta_i$ are i.i.d. variables independent of the initial sample with $E(\zeta_i)=0$ and $E(\zeta_i^2)=E(\zeta_i^3)=1$.
\item[(3)] Obtain a bootstrapped statistic $T_n$ based on the sample $\{ (\tilde Y_i^*,W_i,D_i):i=1,\cdots,n\}$. Repeat this scheme a large number of times, say, $B$ times, the bootstrap critical value at a given level $\alpha$ is  the empirical $(1-\alpha)$-th quantile of the bootstrapped distribution $\{T_{n,j}:j=1,\cdots,B\}$ of the test statistic.
\end{itemize}
Here we use the two point distribution proposed by \cite{mammen1993}:
 $$P(\zeta_i=\frac{1-\sqrt{5}}{2})=\frac{5+\sqrt{5}}{10},
 P(\zeta_i=\frac{1+\sqrt{5}}{2})=\frac{5-\sqrt{5}}{10}.$$

All reported results are based on 1000 simulation runs with 500 bootstrap replications. The sample size $n$ equals $100$ and $200$.
In the following, we report the simulation results at the $\alpha=0.05$ significance level. Without loss of generality, we only consider the case of $X\in R$, i.e. $p=1$ and two dimensions of $Z$: $q=3,~7$. As for nonparametric estimation of $g(x)$, we use the kernel function $K(u)=\frac{3}{4}(1-u^2)$ if $|u|\leq 1$, and $K(u)=0$ otherwise. Thus the order kernel function is $s=2$.

\bigskip
\textbf{Study~1.} Consider the potential outcomes $(Y(1),~Y(0))$ are discrete, taking both low-frequency and high-frequency alternative models into account. The specific data generating processes (DGPs) are as follows:
\begin{eqnarray*}
&& DGP~1:~Y(1)=\mathbb{I}(Y(1)^*>0)~and\ Y(0)=0~with ~Y(1)^* = X+a (\beta^{\top}Z)^3+\epsilon.\\
&& DGP~2:~Y(1)=\mathbb{I}(Y(1)^*>0)~and\ Y(0)=0~with ~Y(1)^* = X+2a \sin(\beta^{\top}Z)+\epsilon.\\
&&The ~observed~ outcome:~Y=DY(1)+(1-D)Y(0)~ with~E(D\mid X)=\frac{\exp(\alpha^{\top}W)}
{1+\exp(\alpha^{\top}W)}.
\end{eqnarray*}

Here  $\beta^{\top}=(\underbrace{1,\cdots,1}_{q})/\sqrt{q}$, $\alpha^{\top}=(1,\underbrace{1,\cdots,1}_{q})/\sqrt{1+q}.$ The observations $\{X_i\}_{i=1}^n\sim i.i.d. N(0,1)$ and $\{Z_i\}_{i=1}^n$, which are independent of $\{X_i\}_{i=1}^n$, are independently generated from
$\mathcal{N}(0,\Sigma),$ $\Sigma=\{\sigma^2_{jj'}\}$,
$\sigma^2_{jj'}=0.5^{|j-j'|}$ for $1\leq j,j'\leq q,~i=1,\cdots,n.$ In the following studies, the errors $\{\epsilon_i\}_{i=1}^n$ are independently drawn from the standard normal distribution. Obviously, the null and alternative hypothesis respectively respond to $a=0$ and $a\neq 0$. Further let $a\in \{0.2,0.4,0.6,0.8,1\}$. Based on DGP~1,  consider the low-frequency alternative model, which is in favour of global smoothing tests. DGP~2 is a high-frequency model under the alternative, which is in favour of local smoothing tests.

Before carrying out simulation procedure, we first  check the sensitivity of  bandwidth selection and choose a reasonable bandwidth. The candidate bandwidths are set to equal  $h_c\sigma n^{-1/4}$ for $h_c\in \{0.6,0.8,\cdots,1.6\}$ and $\sigma=sd(X)$.  To save space, we only investigate the bandwidth impact on $DGP~1$ under $n=200$ and $q=3$. As shown in Table~\ref{band_result}, the different bandwidths have little effect on empirical power, while the empirical size can be still under control when the bandwidth h is not too large. Thus, we choose the bandwidth $h=\sigma n^{-1/4}$ to conduct the following simulation studies and will see that such a choice can also suitable for other models in the numericla studies.
% Table generated by Excel2LaTeX from sheet 'Sheet3'
\begin{table}[htbp]
  \centering
  \caption{The empirical sizes and powers results vary with the bandwidth $h$ and $n=200$, $q=3$ for DGP~1.}
    \begin{tabular}{lrrrrrr}
    \toprule
    $h_c$     & 0.6   & 0.8   & 1     & 1.2   & 1.4   & 1.6 \\
    \midrule
    a=0   & 0.044 & 0.06  & 0.055 & 0.046 & 0.047 & 0.037 \\
    \midrule
    a=1   & 1     & 1     & 1     & 0.999 & 1     & 1 \\
    \bottomrule
        \end{tabular}%
  \label{band_result}%
\end{table}%

Through the simulations results of Study~1 presented in Table \ref{study1}, we have the following observations. First,  the sample sizes reasonably have significant impact on the power performance of $T_n$: large size of sample results in high power, and empirical size of the test is also close to  the significance level. Second, the proposed test  is sensitive to both high and low frequency alternatives.   It is worth mentioning that even for small $a=0.2$, the empirical power of $T_n$ is high enough. To investigate the dimensionality effect of $Z$ on the test performance, we can see that when $q$ increases from $3$ up to $7$, the power performance of $T_n$ is slightly negatively affected, while the empirical size seems to be stable against this dimensionality increasing. Note that for such sample sizes, the total dimension $8$ is already large as $8\approx200^{2/5}=8.33$. Finally, we can see that the test $T_n$ is stil powerful even when dealing with high frequency data. This would be because of the benefit from inheriting its original local smoothing feature.

\begin{table}[htbp]
  \centering
  \caption{Empirical sizes and powers of $T_n$ under study 1.}
      \scalebox{0.95}{

    \begin{tabular}{cccccccccccc}
    \toprule
          & \multicolumn{5}{c}{DGP~1}             &       & \multicolumn{5}{c}{DGP~2} \\
\cmidrule{2-6}\cmidrule{8-12}          & \multicolumn{2}{c}{p+q=4} &       & \multicolumn{2}{c}{p+q=8} &       & \multicolumn{2}{c}{p+q=4} &       & \multicolumn{2}{c}{p+q=8} \\
\cmidrule{2-3}\cmidrule{5-6}\cmidrule{8-9}\cmidrule{11-12}    a     & n=100 & n=200 &       & n=100 & n=200 &       & n=100 & n=200 &       & n=100 & n=200 \\
    \midrule
    0     & 0.048 & 0.046 &       & 0.045 & 0.056 &       & 0.049 & 0.041 &       & 0.06  & 0.049 \\
    0.2   & 0.589 & 0.906 &       & 0.631 & 0.921 &       & 0.243 & 0.336 &       & 0.209 & 0.267 \\
    0.4   & 0.829 & 0.996 &       & 0.844 & 0.992 &       & 0.578 & 0.806 &       & 0.468 & 0.654 \\
    0.6   & 0.929 & 1     &       & 0.887 & 1     &       & 0.849 & 0.938 &       & 0.725 & 0.827 \\
    0.8   & 0.967 & 1     &       & 0.947 & 1     &       & 0.933 & 0.978 &       & 0.82  & 0.886 \\
    1     & 0.985 & 1     &       & 0.952 & 1     &       & 0.946 & 0.979 &       & 0.906 & 0.89 \\
    \bottomrule
    \end{tabular}%
\label{study1}%

  }
\end{table}%

\bigskip
\textbf{Study~2.}
In this study, consider the potential outcomes $(Y(1),~Y(0))$ are continuous. Both low-frequency and high-frequency alternative models are considered via the following data generating processes (DGPs):
\begin{eqnarray*}
&& DGP~3:~Y(1)=2X^2-X+a(\beta^{\top}Z)^3+\epsilon,~and\ Y(0)=0.\\
&& DGP~4:~Y(1)=2X^2-X+ 4a\sin(\beta^{\top}Z)+\epsilon,~and\ Y(0)=0.\\
&&The ~observed~ outcome:~Y=DY(1)+(1-D)Y(0)~ with~E(D\mid X)=\frac{\exp(\alpha^{\top}W)}
{1+\exp(\alpha^{\top}W)}.
\end{eqnarray*}
Here  $\beta^{\top}=(\underbrace{1,\cdots,1}_{q})/\sqrt{q}$, $\alpha^{\top}=(1,\underbrace{1,\cdots,1}_{q})/\sqrt{1+q}.$
 The observations $\{X_i\}_{i=1}^n$, $\{Z_i\}_{i=1}^n$and $\{\epsilon_i\}_{i=1}^n$ are generated as before. The values of $a$ is used to control the deviation from the null hypothesis. DGP~3 and DGP~4 respectively reflect the low-frequency and high-frequency alternative model.

The simulation results are summarized in Table \ref{study2}. Based on the results, we can get similar conclusions as those from Study~1. The main difference between Study 1 and Study 2 is the property of respondse. It seems that $T_n$ performs slightly better when the responds are discrete.

%We can find that the actual finite sample size is close to the nominal size, 0.05.
%
%
%The results show that our test can maintain the significance level well for both dimensions $p=1,q=3$ and $p=1,q=9$.

\begin{table}[htbp]
  \centering
  \caption{Empirical sizes and powers of $T_n$ under study 2.}
      \scalebox{0.95}{

    \begin{tabular}{cccccccccccc}
    \toprule
          & \multicolumn{5}{c}{DGP~3}             &       & \multicolumn{5}{c}{DGP~4} \\
\cmidrule{2-6}\cmidrule{8-12}          & \multicolumn{2}{c}{p+q=4} &       & \multicolumn{2}{c}{p+q=8} &       & \multicolumn{2}{c}{p+q=4} &       & \multicolumn{2}{c}{p+q=8} \\
\cmidrule{2-3}\cmidrule{5-6}\cmidrule{8-9}\cmidrule{11-12}    a   & n=100 & n=200 &       & n=100 & n=200 &       & n=100 & n=200 &       & n=100 & n=200 \\
    \midrule
    0     & 0.048 & 0.047 &       & 0.044 & 0.045 &       & 0.052 & 0.043 &       & 0.052 & 0.052 \\
    0.2   & 0.492 & 0.651 &       & 0.571 & 0.802 &       & 0.317 & 0.395 &       & 0.268 & 0.362 \\
    0.4   & 0.79  & 0.914 &       & 0.808 & 0.925 &       & 0.637 & 0.767 &       & 0.57  & 0.671 \\
    0.6   & 0.86  & 0.96  &       & 0.853 & 0.934 &       & 0.835 & 0.924 &       & 0.743 & 0.846 \\
    0.8   & 0.887 & 0.968 &       & 0.88  & 0.93  &       & 0.931 & 0.969 &       & 0.837 & 0.911 \\
    1     & 0.914 & 0.971 &       & 0.87  & 0.93  &       & 0.956 & 0.974 &       & 0.9   & 0.928 \\
    \bottomrule

    \end{tabular}%

\label{study2}%

  }
\end{table}%

%\section{Real data analysis}
%
% \setcounter{equation}{0}
%
%\newpage

\section{A real data example}

 \setcounter{equation}{0}

 In this section, we consider a data set from AIDS Clinical
 Trials Group Protocol 175 (ACTG175), which can be
 obtained from the R package speff2trial,
 to illustrate the usefulness of the proposed test. This data set
 was popularly analyzed in the literature related to treatment effects,
 such as \cite{hammer1996}, \cite{zhang2008} and \cite{lu2013}.
 There are 2139 HIV-infected subjects
 in ACTG175 that were randomized to four
 different treatment groups with equal probability:
 zidovudine (ZDV) monotherapy, ZDV+didanosine (ddI),
 ZDV+zalcitabine, and ddI monotherapy.
 To get more
 elaborated results, we only consider the subjects under ZDV+didanosine (ddI)
 and ddI monotherapy groups(1083 subjects) in the following analysis.

Hence
 the treatment indicator $D$ is a dummy variable such that
 $D_i=1$ when the $i$-th subject receives ZDV+ddI treatment
 and $D_i=0$ stands for the subject in ddI monotherapy group.
The continuous response $Y$ is CD4 count at
20 $\pm$ 5 weeks post-baseline. Besides, we also consider
 12 baseline covariates acted as $W$, including age, weight, Karnofsky
score, CD4 count at baseline(CD40) and CD8 count at baseline,
hemophilia, homosexual activity,
history of intravenous drug use, race, gender,
antiretroviral history and symptomatic status, which
are also considered in \cite{zhang2008} and
 \cite{lu2013} and more details can be found in their analyses.

The goal of our study is to check whether the
 average treatment effect conditional on whole $W$
 equals  the one conditional on $X$, a subset of
 $W$, i.e.
 \begin{eqnarray}
 H_0:~P(E[Y(1)-Y(0)\mid W]=E[Y(1)-Y(0)\mid X])=1.
 \end{eqnarray}
 \cite{lu2013}
    pointed out that both $age$ and $CD40$ are
   important variables in the contract function $E(Y(1)-Y(0)\mid W)$. We then check whether each individual $age$ or $CD40$ is sufficient or whether they are jointly important. Thus, there are three candidates of $X$ being considered and
  the corresponding null hypotheses are as follows:
 \begin{itemize}
 \item $H_0^1: P(E[Y(1)-Y(0)\mid W]=E[Y(1)-Y(0)\mid X])=1$ with $X=age$;
 \item $H_0^2: P(E[Y(1)-Y(0)\mid W]=E[Y(1)-Y(0)\mid X])=1$ with $X=CD40$;
 \item$H_0^3: P(E[Y(1)-Y(0)\mid W]=E[Y(1)-Y(0)\mid X])=1$ with $X=(age,CD40)$.
 \end{itemize}

 Similarly as in the simulation studies, we use the Epanechnikov
 kernel and then the kernels of  order $s=2p$ are derived from it. The  bandwidth parameter
  $h=\frac{tr(\Sigma_x^{1/2})}{p}n^{-2s}$ with $\Sigma_x=Var(X)$
  and $tr(A)$ being the trace of matrix $A$. Also note that the propensity score $\pi(W)\equiv 0.5$. Hence
   we can get the p-values based on (\ref{tn}) and aforementioned
    wild bootstrap approximation procedure with 500 bootstrap replications.
 The $p$-values of these three tests are respectively $0.06$, $0.082$ and $0.28$
  respectively. Although if we only consider the significance level $0.05$, all three null hypotheses cannot be rejected, the test for $H_{0}^3$ provides  clearer information on the plausibility that a heterogenous
  average treatment effect would be on both $X=(age, CD40)$. Thus, we may not use only either one to model the heterogenous
  average treatment effect to avoid the model misspecification risk.
   This analysis then provides a formal assessment for the model  \cite{lu2013} considered.

\section{Conclusion}

In this study, we consider the testing problem for conditional average treatment effect to explore whether the equivalence relationship between $E(Y(1)-Y(0)\mid W)$ and $E(Y(1)-Y(0)\mid X)$ holds with $X\subset W$. The proposed test has three useful features. That is, it reduces the risk of model misspecification, mitigates the curse of dimensionality through a projection-based construction procedure and achieves the fastest possible rate $n^{-1/2}$ to have the sensitivity to local alternatives. These features make the test  well perform in practice. But this test is not suitable in very large dimension scenarios in the sense that the dimension of covariates is regarded as divergent to infinity. This important research is ongoing.

\renewcommand{\theequation}{A.\arabic{equation}}
\renewcommand{\thesection}{A.\arabic{section}}
\setcounter{section}{0}
\setcounter{equation}{0}
\section*{Appendix}

\section{Preliminary for the proofs.}
\setcounter{equation}{0}

 Before we present the proof, we first define some related quantities.
\begin{itemize}

\item[(1)] $Y_j^*=\frac{D_jY_j}{\pi(X_j,r_0)}-\frac{(1-D_j)Y_j}{1-\pi(X_j,r_0)}=V_j( r_0)Y_j$ with $V_j(r_0)=\frac{D_j}{\pi(X_j,r_0)}-\frac{1-D_j}{1-\pi(X_j, r_0)}$;

\item[(2)] $e_j=Y_j^*-g(X_j)$  and $u_j=Y_j^*-E(Y_j^*\mid W_j)$ with $g(X_j)=E(Y_j^*\mid X_j)$;

\item[(3)] $w_{ij}=\frac{1}{(n-1) }  \mathcal{K}_h\left(X_j-X_i \right)\big/\hat f(X_i)$, $\tilde w_{ij}= \mathcal{K}_h\left(X_j-X_i \right)\big/\hat f(X_i)$ and $\bar w_{ij}=  \mathcal{K}_h\left(X_j-X_i \right)\big/f(X_i)$ with  $\hat f(X_i)=\frac{1}{(n-1)} \sum_{j\neq i}^n  \mathcal{K}_h\left(X_j-X_i \right)$ and $\mathcal{K}_h(\cdot)=\frac{1}{h^p}\mathcal{K}
\left(\frac{\cdot}{h}\right)$;

\item[(4)] $\hat e_j=\hat Y^*_j-\hat g(X_j)$ with $\hat Y^*_j=V_j(\hat r)Y_j$ and $\hat g(X_j)=\sum_{t\neq j} w_{jt}\hat Y^*_t$;

\item[(5)] $\tilde e_j= Y^*_j-\tilde g(X_j)$ with $Y^*_j=V_j(r_0)Y_j$ and $\tilde g(X_j)=\sum_{t\neq j} w_{jt}Y^*_t$;

\item[(6)] $\eta_j=\overline{Y_j}-G_j$ with $\overline{Y_j}=\nabla V_j(r_0)^{\top} Y_j$ and $G_j=G(X_j)=E(\overline Y_j\mid X_j)$. Here $\nabla V_j(r_0)$ stands for the  Partial derivative with respect to $r$;

%\item[(6)] $K_i=K(\frac{X_{1i}
%-x_{10}}{h})$, $\mu_{s_1}(K)=\int u^{s_1}K(u)du,$ $\parallel
%K\parallel_2^2=\int K(u)^2 du$;
%\item[(7)] $Z_i=(X_i,D_i,Y_i)$,
%$f(X_1)$ is the density function of $X_1$, $f_x(X)$ is the density function of $X$, $f_\alpha(\alpha{\top}X)$ is the density function of $\alpha^{\top}X$, and denote $f_j(y(j)\mid x_{1})$ is the point $(Y(j)=y(j), X_1=x_{1})$ in
%the conditional density function of $(Y(j)\mid X_1)$, $j=0,1;$
\item[(7)] $C$ stands for a generic bounded constant and $\Omega$ is the support of $X$.
\end{itemize}

{\textbf{Proof of Theorem~\ref{null_them}.}
Note that under the null hypothesis $e_i=u_i$ such that $E(u_i\mid W_i)=0$. Recall that $\hat e_j=\hat Y^*_j-\hat g(X_j)$ and $\tilde e_j= Y^*_j-\tilde g(X_j)$. Hence $T_n$ can be decomposed as
\begin{eqnarray*}
&&T_n=\frac{1}{n(n-1)}\sum_{i=1}^n\sum_{j\neq i}\hat e_i\hat e_j \mathcal{B}_{ij}\nonumber\\
&&=\frac{1}{n(n-1)}\sum_{i=1}^n\sum_{j\neq i} u_i u_j\mathcal{B}_{ij} +\frac{2}{n(n-1)}\sum_{i=1}^n\sum_{j\neq i} u_i(\tilde e_j-u_j) \mathcal{B}_{ij} +\frac{1}{n(n-1)}\sum_{i=1}^n\sum_{j\neq i}  (\tilde e_i-u_i)(\tilde e_j- u_j)\mathcal{B}_{ij}\\
&&+\frac{2}{n(n-1)}\sum_{i=1}^n\sum_{j\neq i} u_i(\hat e_j-\tilde e_j) \mathcal{B}_{ij} +\frac{2}{n(n-1)}\sum_{i=1}^n\sum_{j\neq i}  (\hat  e_i-\tilde e_i)(\tilde e_j- u_j)\mathcal{B}_{ij}\\
 &&+\frac{1}{n(n-1)}\sum_{i=1}^n\sum_{j\neq i}  (\hat e_i -\tilde e_i)(\hat e_j-\tilde e_j)\mathcal{B}_{ij}\\
&&:=\frac{1}{n(n-1)}\sum_{i=1}^n\sum_{j\neq i} u_i u_j\mathcal{B}_{ij} +I_{n1}+\cdots+I_{n5}.
\end{eqnarray*}

We then deal with $I_{ni}$ for $1\le i\le 5$ in Propositions~$1\sim5$ separately. The proofs of these propositions will be given later.

{\textbf{Proposition~1.} $I_{n1}=-\frac{2}{n(n-1)}
\sum_{i=1}^n\sum_{j\neq i}u_iu_j\tilde{\mathcal{B}}_{ij}-\frac{1}{n}E[(u_1^2\bar w_{12}+u_2^2\bar w_{21})
\mathcal{B}_{12}]+o_p\left(\frac{1}{n}\right)$ with $\tilde{\mathcal{B}}_{ij}=\frac{1}{2}
 E\left(\bar w_{ti}
\mathcal{B}_{jt}+\bar w_{tj}\mathcal{B}_{it}\mid W_i,W_j\right).$}

{\textbf{Proposition~2.} $I_{n2}=\frac{1}{n(n-1)}
\sum_{i=1}^n\sum_{j\neq i}u_iu_j\bar{\mathcal{B}}_{ij}+\frac{1}{3n}E(u_1^2\bar w_{21}
\bar w_{31}\mathcal{B}_{23}
+u_2^2\bar w_{12}\bar w_{32}\mathcal{B}_{13}+u_3^2\bar w_{12}\bar w_{23}\mathcal{B}_{12})+o_p\left(\frac{1}{n}\right)$
 with $\bar{\mathcal{B}}_{ij}=\frac{1}{2}
 E\left(\bar w_{ti}\bar w_{kj}\mathcal{B}_{tk}
 +\bar w_{ki}\bar w_{tj}\mathcal{B}_{kt}\mid W_i,W_j\right)$.}

{\textbf{Proposition~3.} $I_{n3}=(\hat r-r_0)^{\top}\frac{2}{n}
\sum_{i=1}^nu_iH_3(W_i)+o_p\left(\frac{1}{n}\right)$ with $H_3(W_i)=E(\eta_j{\mathcal{B}}_{ij}\mid W_i)$.}

{\textbf{Proposition~4.} $I_{n4}=-(\hat r-r_0)^{\top}\frac{2}{n }
\sum_{i=1}^n u_iH_4(W_i)+o_p\left(\frac{1}{n}\right)$ with
 $H_4(W_i)=\frac{1}{2}E\left[(\eta_j\bar w_{ti}
+\eta_t\bar w_{ji})\mathcal{B}_{tj}\mid W_i\right]$.

{\textbf{Proposition~5.} $I_{n5}=(\hat r-r_0)^{\top}A(\hat r -r_0)+o_p\left(\frac{1}{n}\right)$ with $A=E(\eta_1\eta_2^{\top}{\mathcal{B}}_{12})$.}

Based on these propositions, by the results about the standard first-order
 degenerate U-statistic in \cite{serfling1980}, the proof is then finished.
  \begin{flushright}
 $\Box$
 \end{flushright}

\bigskip
\bigskip
\bigskip
{\textbf{Proof of Theorem~\ref{alter_them}.}
Under the alternative hypothesis $H_{n1}$, $T_n$ can also be decomposed as
\begin{eqnarray*}
&&T_n=\frac{1}{n(n-1)}\sum_{i=1}^n\sum_{j\neq i}\hat e_i\hat e_j \mathcal{B}_{ij}\nonumber\\
&&=\frac{1}{n(n-1)}\sum_{i=1}^n\sum_{j\neq i} e_i e_j\mathcal{B}_{ij} +\frac{2}{n(n-1)}\sum_{i=1}^n\sum_{j\neq i} e_i(\tilde e_j-e_j) \mathcal{B}_{ij} +\frac{1}{n(n-1)}\sum_{i=1}^n\sum_{j\neq i}  (\tilde e_i-e_i)(\tilde e_j- e_j)\mathcal{B}_{ij}\\
&&+\frac{2}{n(n-1)}\sum_{i=1}^n\sum_{j\neq i} e_i(\hat e_j-\tilde e_j) \mathcal{B}_{ij} +\frac{2}{n(n-1)}\sum_{i=1}^n\sum_{j\neq i}  (\hat  e_i-\tilde e_i)(\tilde e_j- e_j)\mathcal{B}_{ij} \\
&&+\frac{1}{n(n-1)}\sum_{i=1}^n\sum_{j\neq i}  (\hat e_i -\tilde e_i)(\hat e_j-\tilde e_j)\mathcal{B}_{ij}\\
&&:=J_{n0}+J_{n1}+\cdots+J_{n5}.
\end{eqnarray*}

 Similarly as those for $J_{ni}$ for $1\le i\le 5$, we will prove the following Propositions $6\sim11$ under Assumptions $1\sim9$.

{\textbf{Proposition~6.} $J_{n0}= \frac{1}{n(n-1)}
\sum_{i=1}^n\sum_{j\neq i}u_iu_j\tilde{\mathcal{B}}_{ij}+
2a_n\frac{1}{n}\sum_{i=1}^n  u_i \tilde H_{00}(W_i) +a_n^2\mu_0$ with $H_{00}(W_i)=E(H_j\mathcal{B}_{ij\mid W_i})$
 and $\mu_0=E(H_1H_2\mathcal{B}_{12})$.

{\textbf{Proposition~7.} $J_{n1}=-\frac{2}{n(n-1)}
\sum_{i=1}^n\sum_{j\neq i}u_iu_j\tilde{\mathcal{B}}_{ij}
-\frac{1}{n}E[(u_1^2\bar w_{12}+u_2^2\bar w_{21})
\mathcal{B}_{12}]+o_p\left(\frac{1}{n}\right)
-2a_n\frac{1}{n}\sum_{i=1}^nu_iH_{10}(W_i)+o_p
\left(a_n^2\right)$ with $H_{10}(W_i)=\frac{1}{2}E\left[(H_j\bar w_{ti}+H_t\bar w_{ji})
\mathcal{B}_{tj}\mid W_i\right]$.

{\textbf{Proposition~8.} $J_{n2}= \frac{1}{n(n-1)}
\sum_{i=1}^n\sum_{j\neq i}u_iu_j\bar{\mathcal{B}}_{ij}
+\frac{1}{3n}E(u_1^2\bar w_{21}
\bar w_{31}\mathcal{B}_{23}
+u_2^2\bar w_{12}\bar w_{32}\mathcal{B}_{13}+u_3^2\bar w_{12}\bar w_{23}\mathcal{B}_{12})
+o_p\left(\frac{1}{n}\right) +o_p\left(\frac{a_n}{\sqrt n}\right)+o_p(a_n^2)$.}

{\textbf{Proposition~9.} $J_{n3}= (\hat r-r_0)^{\top}\frac{2}{n}
\sum_{i=1}^nu_iH_3(W_i)+a_n(\hat r-r_0)^{\top}\alpha+o_p\left(\frac{a_n}{ \sqrt n}\right)$ with $\alpha=E(H_1\eta_2\mathcal{B}_{12})$.

{\textbf{Proposition~10.} $J_{n4}=-(\hat r-r_0)^{\top}\frac{2}{n }
\sum_{i=1}^n u_iH_4(W_i)+o_p\left(\frac{1}{n}\right)+o_p\left(\frac{1}{ n}\right)+o_p\left(\frac{a_n}{\sqrt n}\right)$.

 {\textbf{Proposition~11.} $J_{n5}=(\hat r-r_0)^{\top}A(\hat r -r_0)+o_p\left(\frac{1}{n}\right)$.}

Based on these propositions, the conclusion in Theorem~1 and the central limit theorem, we can prove Theorem~2.
  \begin{flushright}
 $\Box$
 \end{flushright}

 \bigskip

{\textbf{Proof of Proposition~1.}}
Recall that under the null hypothesis, $Y_j^*=g(X_j)+u_j$ with $E(u_j\mid W_j)=0$, it follows that
\begin{eqnarray}\label{e_dec}
\tilde e_j- u_j=g(X_j)-\sum_{t\neq j}w_{jt}g(X_t)-\sum_{t\neq j}w_{jt}u_t.
\end{eqnarray}
Thus we can decompose $I_{n1}$ as
\begin{eqnarray*}
&&I_{n1}=\frac{2}{n(n-1)}\sum_{i=1}^n\sum_{j\neq i}  u_i (\tilde e_j-u_j)\mathcal{B}_{ij}\\
&&\qquad=\frac{2}{n(n-1)}\sum_{i=1}^n\sum_{j\neq i}  u_i \Big (g(X_j)-\sum_{t\neq j}w_{jt}g(X_t)\Big )\mathcal{B}_{ij}-\frac{2}{n(n-1)}
\sum_{i=1}^n\sum_{j\neq i}\sum_{t\neq j}  u_iw_{jt}u_t\mathcal{B}_{ij}\\
&&\qquad:=2I_{n11}-2I_{n12}.
\end{eqnarray*}

Consider the term $I_{n11}$ first. Without loss of generality, here we consider $s=2$ for ease of exposition. Since $\mathcal{K}$ is a kernel of order $s=2$,
by standard nonparametric theory in literature, e.g. \cite{hardle2012}, the following formula holds uniformly:
\begin{eqnarray*}
g(X_j)-\sum_{t\neq j}w_{jt}g(X_t)=h^2\mathcal{K}_2\sum_{l=1}^p\frac{2\frac{
\partial g(X_j)}{\partial X_{jl}}\frac{
\partial f(X_j)}{\partial X_{jl}}+f(X_j) \frac{
\partial^2 g(X_j)}{\partial X_{jl}\partial X_{jl}}}{2f(X_j)}+o_p(h^2):=h^2F_j+o_p(h^2).
\end{eqnarray*}
Here $\mathcal{K}_2I_p= \int uu^{\top}\mathcal{K}(u)du$ with {$I_p$ being a $p\times p$ identity matrix.}
Further by $E(u_i\mid W_i)=0$, it follows that
\begin{eqnarray*}
&&E(I_{n11}^2)=\frac{1}{n^2(n-1)^2}\sum_{i=1}^n\sum_{j\neq i}  \sum_{i'=1}^n\sum_{j'\neq i'}E\bigg[
u_iu_{i'} \{g(X_j)-\sum_{t\neq j}w_{jt}g(X_t)\}\\
&&\qquad\qquad\qquad\qquad\qquad\qquad\qquad\qquad\qquad\qquad\times\{g(X_{j'})-\sum_{t'\neq j'}w_{j't'}g(X_{t'})\}\mathcal{B}_{ij}\mathcal{B}_{i'j'}\bigg]\\
&&\qquad=h^{4}\frac{1}{n^2(n-1)^2}\sum_{i=1}^n
\sum_{j\neq i}  \sum_{j'\neq i}
E\left[u_i^2  F_jF_{j'} \mathcal{B}_{ij}\mathcal{B}_{i'j'}\right]+R_n =O\left(\frac{h^{4}}{n}\right)+o\left(\frac{h^{4}}{n}\right).
\end{eqnarray*}
Thus we have $I_{n11}=O_p(\frac{h^2}{\sqrt{n}})$. By the assumption $nh^{4}\rightarrow 0$,  $I_{n11}=o_p\left(\frac{1}{n}\right)$.
When $s>2$, we can similarly get $I_{n11}=O_p(\frac{h^{s}}{\sqrt{n}})=o_p\left(\frac{1}{n}\right)$.

As for $I_{n12}$, we can rewrite it as
\begin{eqnarray*}
&&I_{n12}= \frac{1}{n(n-1)^2}
\sum_{i=1}^n\sum_{j\neq i} u_i^2\tilde w_{ji} \mathcal{B}_{ij}+\frac{1}{n(n-1)}
\sum_{i=1}^n\sum_{j\neq i}\sum_{t\neq j,i}  u_iw_{jt}u_t\mathcal{B}_{ij}\\
&&\qquad:=I_{n121}+I_{n122}.
\end{eqnarray*}

%By the Lemma A.1 in \cite{abrevaya2015}, given Assumptions $4,6,7$, $w_{ij}$ is asymptotic symmetric, i.e.
%\begin{eqnarray}\label{sys_w}
%|w_{ij}-w_{ji}|\leq \frac{C_n}{nh^k}|\mathcal{K}\left(\frac{X_i-X_j}{h}\right)|.
%\end{eqnarray}
Given Assumptions $5$ and $6$, $E|u_i^2w_{ij}\mathcal{B}_{ij}|\leq E|u_i^2|<\infty$ and noting that
\begin{eqnarray}\label{fest}
\frac{1}{\hat f(X_i)}=\frac{1}{f(X_i)}+O_p(\hat f(X_i)-f(X_i))=\frac{1}{f(X_i)}+o_p(1),
 \end{eqnarray}
employing the standard U-statistic theory and the results about nonparametric estimation, we have $(n-1)I_{n121} -\frac{1}{2}E[(u_1^2\bar w_{12}+u_2^2\bar w_{21})\mathcal{B}_{12}]=o_p(1)$ with $\bar w_{12}=\frac{\mathcal{K}_h(X_1-X_2)}{f(X_1)}.$ Thus $I_{n121}=\frac{1}{2n}E[(u_1^2\bar w_{12}+u_2^2\bar w_{21})\mathcal{B}_{12}]+o_p\left(\frac{1}{n}\right)$.

Consider $I_{n122}$. Write it as  $I_{n122}=\frac{n-2}{n-1}U_n$ where
\begin{eqnarray*}
&&U_n=\frac{1}{C_n^3}\mathop{\sum\sum\sum}_{1\leq i<j<t\leq n}I^s(\chi_i,\chi_j,\chi_t).
\end{eqnarray*}
Here $U_n$ is an  U-statistic of order 3 and $\chi_i=(W_i,u_i)$, $I^s(\chi_i,\chi_j,\chi_t)=(I_{ijt}+I_{itj}+I_{jit}+
I_{jti}+I_{tij}+I_{tji})/6$
is the kernel with $I_{ijt}=u_i\tilde w_{jt}u_t\mathcal{B}_{ij}$ and $\tilde w_{jt}=\frac{\mathcal{K}_h(X_j-X_t)}{\hat f(X_j)}$.

Note that $E(u_i\mid W_1,\cdots,W_n)=0$. It follows that $E(I^s(\chi_i,\chi_j,\chi_t)\mid \chi_i)=0$ and
$E(I^s(\chi_i,\chi_j,\chi_t)\mid \chi_i,\chi_j)=\frac{1}{6}u_iu_jE\left(\tilde w_{ti}
\mathcal{B}_{jt}+\tilde w_{tj}\mathcal{B}_{it}\mid W_i,W_j\right)\neq 0$. Thus $U_n$ is a  degenerate U-statistic. Given assumptions~4 and 5, we can also observe that $E(u_i^2u_t^2\tilde w_{jt}^2\mathcal{B}_{ij}^2)\leq CE(u_i^2u_t^2)<\infty$. Therefore, by using the results collected in \cite{serfling1980} and (\ref{fest}), we have
\begin{eqnarray*}
&&U_n=\frac{2}{n(n-1)}\mathop{\sum\sum}_{1\leq i<j\leq n}u_iu_j\tilde{\mathcal{B}}_{ij}+o_p\left(
\frac{1}{n}\right)=\frac{1}{n(n-1)}{\sum_{i=1}^n\sum_{j\neq i}}u_iu_j\tilde{\mathcal{B}}_{ij}+o_p\left(
\frac{1}{n}\right)
\end{eqnarray*}
 Here $\tilde{\mathcal{B}}_{ij}=\frac{1}{2}
 E\left(\bar w_{ti}
\mathcal{B}_{jt}+\bar w_{tj}\mathcal{B}_{it}\mid W_i,W_j\right).$
Taking all the asymptotic results about the terms in the decomposition of $I_{n1}$ into account, proposition~1 is proved.
 \begin{flushright}
 $\Box$
 \end{flushright}

{\textbf{Proof of Proposition~2.}}
Based on (\ref{e_dec}), $I_{n2}$ can be decomposed as
\begin{eqnarray*}
&&I_{n2}=\frac{1}{n(n-1)}\sum_{i=1}^n\sum_{j\neq i}  (\tilde e_i-u_i)(\tilde e_j- u_j)\mathcal{B}_{ij} \\
&&\qquad=\frac{1}{n(n-1)}\sum_{i=1}^n\sum_{j\neq i}  \{g(X_i)-\sum_{t\neq i}w_{it}g(X_t)\}\{g(X_j)-\sum_{k\neq j}w_{jk}g(X_k)\}\mathcal{B}_{ij} \\
&&\qquad-2\frac{1}{n(n-1)}\sum_{i=1}^n\sum_{j\neq i} \sum_{k\neq j}\{g(X_i)-\sum_{t\neq j}w_{it}g(X_t)\}w_{jk}u_k\mathcal{B}_{ij} \\
&&\ \qquad+\frac{1}{n(n-1)}\sum_{i=1}^n\sum_{j\neq i}  \sum_{t\neq i}\sum_{k\neq j}w_{it}u_tw_{jk}u_k\mathcal{B}_{ij} \\
&&\qquad:=I_{n21}-2I_{n22}+I_{n23}.
\end{eqnarray*}
As proved above, with $g(X_j)-\sum_{t\neq j}w_{jt}g(X_t)=O_p(h^s)$, we can get $I_{n21}=O_p(h^{2s}).$ Further, similarly as the proof for $I_{n11}$, we can show $I_{n22}=O_p(\frac{h^s}{\sqrt{n}})=
o_p\left(\frac{1}{n}\right).$

For the term $I_{n23}$, we also have the decomposition as
\begin{eqnarray*}
&&I_{n23}=\frac{1}{n(n-1)^3}\sum_{i=1}^n\sum_{j\neq i}  \sum_{t\neq i,j} \tilde w_{it}\tilde w_{jt}u_t^2\mathcal{B}_{ij}+\frac{1}{n(n-1)^3}\sum_{i=1}^n\sum_{j\neq i}  \sum_{t\neq i,j}\sum_{k\neq i,j,t}\tilde w_{it}\tilde w_{jk}u_tu_k\mathcal{B}_{ij}+R_n\\
&&\qquad:=I_{n231}+I_{n232}+R_n.
\end{eqnarray*}
Since $E(u_i\mid W_1,\cdots,W_n)=0$, the leading term of $I_{n23}$ is $I_{n231}+I_{n232}$ and $R_n$ is a higher order term with $E(R_n)=0$. In the following, we focus on investigating the asymptotic behaviour of $I_{n231}+I_{n232}$.

Rewrite $I_{n231}$ as a U-statistic of order 3: $I_{n231}=\frac{n-2}{(n-1)^2}V_n$ with
\begin{eqnarray*}
&&V_n=\frac{1}{C_n^3}\mathop{\sum\sum\sum}_{1\leq i<j<t\leq n}V^s(\chi_i,\chi_j,\chi_t).
\end{eqnarray*}
Here $V_n$ is an U-statistic of order 3 and $\chi_i=(W_i,u_i)$, $V^s(\chi_i,\chi_j,\chi_t)=(V_{ijt}+V_{itj}+
I_{jti}  )/3$
is the kernel with $V_{ijt}= \tilde w_{it}\tilde w_{jt}u_t^2\mathcal{B}_{ij}$.

Notably, $E|\tilde w_{it}\tilde w_{jt}u_t^2\mathcal{B}_{ij}|\leq CE(u_t^2)<\infty$, and thus employing the U-statistic theory in \cite{serfling1980}, we can obtain that $V_n \stackrel{\mathcal{P}}{\longrightarrow}\mu_v$ with $\mu_v=E(V^s(\chi_i,\chi_j,\chi_t))$. That implies $I_{n231}=\frac{1}{n}\mu_v+o_p\left(
\frac{1}{n}\right)$.

For $I_{n232},$ we can also write it as $I_{n232}=\frac{(n-2)(n-3)}{(n-1)^2}V_n^*$ with
\begin{eqnarray*}
&&V^*_n=\frac{1}{C_n^4}\mathop{\sum\sum\sum\sum}_{1\leq i<j<t<k\leq n}V^{*s}(\chi_i,\chi_j,\chi_t,\chi_k).
\end{eqnarray*}
Here $V_n^*$ is an  U-statistic of order 4  and $V^{*s}(\chi_i,\chi_j,\chi_t,\chi_k)=\frac{1}{24}\sum_{(a)}
V^*(\chi_{i1},\chi_{i2},\chi_{i3},\chi_{i4})$ is a symmetric kernel where $\sum_{a}$ denotes summation over the $4!$ permutations $(i_1,\cdots,i_4)$ of $(i,j,t,k)$.

Given assumptions~4 and 5, on the analogy of investigating $I_{n122}$, we can obtain $V^*_n$ is a  degenerate U-statistic and
\begin{eqnarray*}
&&V^*_n=\frac{2}{n(n-1)}\mathop{\sum\sum}_{1\leq i<j\leq n}u_iu_j\bar{\mathcal{B}}_{ij}+o_p\left(
\frac{1}{n}\right)=\frac{1}{n(n-1)}{\sum_{i=1}^n\sum_{j\neq i}}u_iu_j\bar{\mathcal{B}}_{ij}+o_p\left(
\frac{1}{n}\right)
\end{eqnarray*}
 Here $\bar{\mathcal{B}}_{ij}=\frac{1}{2}
 E\left(\bar w_{ti}\bar w_{kj}\mathcal{B}_{tk}
 +\bar w_{ki}\bar w_{tj}\mathcal{B}_{kt}\mid W_i,W_j\right).$
Altogether,  we can conclude the proof of Proposition~2.
 \begin{flushright}
 $\Box$
 \end{flushright}

{\textbf{Proof of Proposition~3.}}
Note that $\hat e_j-\tilde e_j=(V_j(\hat r)-V_j(r_0))Y_j-\sum_{t\neq j}w_{jt}(V_t(\hat r)-V_t(r_0))Y_t.$
%=[\nabla V_j(r)Y_j-\sum_{t\neq j}w_{jt}\nabla V_t(r)Y_t]^{\top}(\hat r -r)+ (\hat r-r)o_p(1)$.
 Let $\overline{Y_j}=\nabla V_j(r_0)^{\top} Y_j$, $G_j=G(X_j)=E(\overline Y_j\mid X_j)$,  $\eta_j=\overline{Y_j}-G_j$ with $E(\eta_i\mid W_i)\neq 0$.

By the Taylor  expansion around $r_0$ and $\hat r-r_0=O_p(n^{-1/2})$,
it follows that
\begin{eqnarray*}
&&I_{n3}=(\hat r-r_0)^{\top}\frac{2}{n(n-1)}\sum_{i=1}^n\sum_{j\neq i}u_i (\overline Y_j-\sum_{t\neq j}w_{jt}\overline Y_t)   \mathcal{B}_{ij}+R_n\\
&&= (\hat r-r_0)^{\top} \frac{2}{n(n-1)}\sum_{i=1}^n \sum_{j\neq i} u_i (G_j-\sum_{t\neq j}w_{jt}G_t) \mathcal{B}_{ij}+(\hat r-r_0)^{\top}\frac{2}{n(n-1)}\sum_{i=1}^n \sum_{j\neq i}  u_i\eta_j \mathcal{B}_{ij}\\
&&- (\hat r-r_0)^{\top}\frac{2}{n(n-1)}\sum_{i=1}^n \sum_{j\neq i}\sum_{t\neq j}  u_i w_{jt}\eta_t \mathcal{B}_{ij} +R_n\\
&&:=2(\hat r-r_0)^{\top}I_{n31}+2(\hat r-r_0)^{\top}I_{n32}-2(\hat r-r_0)^{\top}I_{n33}+R_n.
\end{eqnarray*}
Here $R_n=o_p(I_{n31}+I_{n32}-I_{n33})$ is a higher order term.

For the term $I_{n31}$, we can show that
\begin{eqnarray*}
E(I_{n31}^2)=\frac{1}{n^2(n-1)^2}\sum_{i=1}^n \sum_{j\neq i}\sum_{i'=1}^n \sum_{j'\neq i'} E( u_i u_i' (G_j-\sum_{t\neq j}w_{jt}G_t)   (G_{j'}-\sum_{t'\neq j'}w_{j't'}G_{t'})^{\top} \mathcal{B}_{ij}\mathcal{B}_{i'j'}).
\end{eqnarray*}
Noting that $\{W_i\}_{i=1}^n$ are i.i.d. samples, $E(u_i\mid W_i)=0$, thus only the terms $i=i'$ have non-zero expectations, therefore,
\begin{eqnarray*}
E(I_{n31}^2)=\frac{1}{n^2(n-1)^2}\sum_{i=1}^n\sum_{j\neq i} \sum_{j'\neq i} E( u_i^2 (G_j-\sum_{t\neq j}w_{jt}G_t)  (G_{j'}-\sum_{t'\neq j'}w_{j't'}G_{t'})\mathcal{B}_{ij}\mathcal{B}_{i'j'}).
\end{eqnarray*}
Further, note that $G_j-\sum_{t\neq j}w_{jt}G_t$ is also a bias term. Hence similarly as the proof for  $I_{n11}$, we have $E(I_{n31}^2)=O(\frac{h^{2s}}{n})
=o_p\left(\frac{1}{n}\right).$
That implies $(\hat r-r_0)^{\top}I_{n31}=o_p\left(\frac{1}{n}\right)$.

As for $I_{n33}$, we can rewrite it as
\begin{eqnarray*}
&&I_{n33}=\frac{1}{n(n-1)^2}\sum_{i=1}^n \sum_{j\neq i}   u_i \eta_i \tilde w_{ji}\mathcal{B}_{ij}+\frac{1}{n(n-1)^2}\sum_{i=1}^n \sum_{j\neq i}\sum_{t\neq i,j}  u_i \tilde w_{jt}\eta_t \mathcal{B}_{ij}\\
&&:=I_{n331}+I_{n332}.
\end{eqnarray*}
Similarly as the discussion on $I_{n231}$,  we can show that $I_{n331}= \frac{1}{n}\mu_{331}+o_p\left(\frac{1}{n}\right)$ with $\mu_{331}=\frac{1}{2}E(u_i\eta_i\bar w_{ji}\mathcal{B}_{ij}+u_j\eta_j\bar w_{ij}\mathcal{B}_{ij})$ and then $(\hat r-r_0)^{\top}I_{n331}=o_p\left(\frac{1}{n}\right).$
It is easy to see that
\begin{eqnarray*}
E(I_{n332}^2)=\frac{1}{n^2(n-1)^4}\sum_{i=1}^n
\sum_{j\neq i}\sum_{t\neq i,j}\sum_{i'=1}^n \sum_{j'\neq i'}\sum_{t'\neq i', j'} E( u_iu_{i'} \tilde w_{jt} \tilde w_{j't'}\eta_t \eta_{t'}\mathcal{B}_{ij}  \mathcal{B}_{i'j'} ).
\end{eqnarray*}
Since $E(\bar w_{jt}\eta_t\mid W_1,\cdots,W_n)=0$, only the terms $i=i',~t=t'$ have, as discussed before, non-zero expectation, thus
\begin{eqnarray*}
E(I_{n332}^2)=\frac{1}{n^2(n-1)^4}\sum_{i=1}^n
\sum_{j\neq i}\sum_{t\neq i,j} \sum_{j'\neq i}  E( u_i^2 \tilde w_{jt} \tilde w_{j't}\eta_t^2 \mathcal{B}_{ij}  \mathcal{B}_{ij'} )
\end{eqnarray*}
Hence we can obtain $E(I_{n332}^2)=O(\frac{n^4}{n^6})=O(\frac{1}{n^{2}})$. Combining that $\hat r-r=O_p(n^{-1/2})$, $(\hat r-r_0)^{\top}I_{n332}=o_p\left(\frac{1}{n}\right).$

As for $I_{n32}$, we can also rewrite it in the form of U-statistic:
\begin{eqnarray}
I_{n32}=\frac{2}{n(n-1)}\mathop{\sum}_{1\leq i<j\leq n}h_3(\chi_i,\chi_j).
\end{eqnarray}
Here $h_{3}(\chi_i,\chi_j)=\frac{1}{2}(u_i\eta_j+u_j\eta_i)\mathcal{B}_{ij}$ with $\chi_i=(W_i,u_i,D_i)$. Note that $E(h_3(\chi_i,\chi_j)\mid \chi_i)=\frac{1}{2}u_iH_3(W_i)\neq 0$ with $H_3(W_i)=E(\eta_j\mathcal{B}_{ij}\mid W_i)$.
Obviously, $\xi_1=E(u_1^2H_3(W_1)^2)>0$, thus $I_{n32}$ is not degenerate. Further, we have $E(h_3(\chi_i,\chi_j)^2)\leq CE(u_1^2\eta_1^2)$. Then by Holder's inequality and assumption~5, we have $E(u_1^2\eta_2^2)\leq (E(u_1^4))^{1/2}(E(\eta_2^4))^{1/2}< \infty$. Again U-statistic theory in \cite{serfling1980} leads to an asymptotically linear representation of $I_{n32}$ as
\begin{eqnarray}
I_{n32}=\frac{2}{n}
\sum_{i=1}^nu_iH_3(W_i)+o_p\left(\frac{1}{\sqrt n}\right).
\end{eqnarray}
That implies the higher order term $R_n=o_p\left(\frac{1}{n}\right)$.
Altogether, Proposition~3 is proved.  \begin{flushright}
 $\Box$
\end{flushright}

{\textbf{Proof of Proposition~4.}}
As for $I_{n4}$, we can decompose it as
\begin{eqnarray*}
&&I_{n4}=\frac{2}{n(n-1)}\sum_{i=1}^n\sum_{j\neq i}  (\hat e_i-\tilde e_i)\{g(X_j)-\sum_{t\neq j}w_{jt}Y_t^*\}\mathcal{B}_{ij}\\
&&\qquad=(\hat r-r_0)^{\top}\frac{2}{n(n-1)}\sum_{i=1}^n\sum_{j\neq i} \{\overline Y_i-\sum_{k\neq i}w_{it}\overline Y_k\} \{g(X_j)-\sum_{t\neq j}w_{jt}Y_t^*\}\mathcal{B}_{ij}+R_n.
\end{eqnarray*}
Here $R_n$ is a higher order term.
Further, we can decompose $\overline Y_i-\sum_{k\neq i}w_{ik}\overline Y_t=\eta_i+G_i-\sum_{k\neq i}w_{ik}\bar Y_k$, then $I_{n4}$ can be rewritten as
 \begin{eqnarray*}
&&I_{n4}=(\hat r-r_0)^{\top}\frac{2}{n(n-1)}\sum_{i=1}^n\sum_{j\neq i} \{ G_i-\sum_{k\neq i}w_{it}\overline Y_k\} \{g(X_j)-\sum_{t\neq j}w_{jt}Y_t^*\}\mathcal{B}_{ij}\\
&&\qquad\qquad\qquad+(\hat r-r_0)^{\top}\frac{2}{n(n-1)}\sum_{i=1}^n\sum_{j\neq i} \eta_i \{g(X_j)-\sum_{t\neq j}w_{jt}Y_t^*\}\mathcal{B}_{ij}+R_n\\
&&\qquad :=2(\hat r-r_0)^{\top}I_{n41}+2(\hat r-r_0)^{\top}I_{n42}+R_n.
\end{eqnarray*}

By the  standard results in kernel estimation, see e.g. \cite{abrevaya2015}, we have
\begin{eqnarray}\label{unif_non}
&&\mathop{\sup}_{X_i\in \Omega}|G_i-\sum_{k\neq i}w_{ik}\overline Y_k|=O_p\left(h^s+\sqrt{\frac{\log n}{nh^p}}\right),\\
&&\mathop{\sup}_{X_j\in \Omega}|g(X_j)-\sum_{t\neq j}w_{jt}Y_t^*|=O_p\left(h^s+\sqrt{\frac{\log n}{nh^p}}\right).
\end{eqnarray}

Thus we can derive that
\begin{eqnarray*}
|I_{n41}|\leq \mathop{\sup}_{X_i\in \Omega}|G_i-\sum_{k\neq i}w_{it}\overline Y_k|\mathop{\sup}_{X_j\in \Omega}|g(X_j)-\sum_{t\neq j}w_{jt}Y_t^*|\sum_{i=1}^n\sum_{j\neq i} |\mathcal{B}_{ij}|
\end{eqnarray*}
Since $|\mathcal{B}_{ij}|\leq 1$, then under assumption~7, it follows that $I_{n41}=o_p\left(\frac{1}{\sqrt{n}}\right)$.
So that $(\hat r-r_0)^{\top}I_{n41}=o_p\left(\frac{1}{n}\right)$.

As for $I_{n42}$, note that $Y_t^*=g(X_t)+u_t$  and then
 \begin{eqnarray*}
&&I_{n42}=\frac{1}{n(n-1)}\sum_{i=1}^n\sum_{j\neq i} \eta_i \{g(X_j)-\sum_{t\neq j}w_{jt}g(X_t)\}\mathcal{B}_{ij}-\frac{1}{n(n-1)^2}\sum_{i=1}^n\sum_{j\neq i}  \sum_{t\neq j,i}\eta_i \tilde w_{jt}u_t\mathcal{B}_{ij}\\
&&\qquad-\frac{1}{n(n-1)^2}\sum_{i=1}^n\sum_{j\neq i}  \eta_iu_i \tilde w_{ji}\mathcal{B}_{ij}\\
&&:=I_{n421}-I_{n422}-I_{n423}.
 \end{eqnarray*}
As the bias term $\{g(X_j)-\sum_{t\neq j}w_{jt}g(X_t)\}=O_p(h^s)$, under assumption~7, we can easily obtain that $I_{n421}=O_p(h^s)=o_p\left(\frac{1}
{\sqrt{n}}\right)$. Thus we have $(\hat r-r_0)I_{n421}=o_p\left(\frac{1}
{n}\right)$. Further, $I_{n423}=O_p\left(\frac{1}
{n}\right)=o_p\left(\frac{1}
{\sqrt n}\right)$ is obviously derived. Thus we have $(\hat r-r_0)I_{n423}=o_p\left(\frac{1}
{n}\right)$ as well.

 Next, we deal with $I_{n422}$. As the arguments are very similar, we then only give an outline. Rewrite it as  an $U$-statistic
 $I_{n422}=\frac{n-2}{n-1}V^{\star}_n$ with
 \begin{eqnarray}
 V_n^{\star}=\frac{1}{C_{n}^3}\mathop{\sum\sum\sum}_{1\leq i<j<t\leq n}V^{\star s}(\tilde \chi_i,\tilde \chi_j,\tilde \chi_t)
 \end{eqnarray}
Here $U_n$ is an U-statistic of order 3 and $\tilde \chi_i=(W_i,u_i,\eta_i)$, $V^{\star s}(\tilde\chi_i,\tilde\chi_j,\tilde\chi_t)=(V^\star_{ijt}+V^\star_{itj}
+V^\star_{jit}+V^\star_{jti}+V^\star_{tij}+V^\star_{tji})/6$
is the kernel with $V^{\star}_{ijt}= \eta_i \bar w_{jt}u_t\mathcal{B}_{ij}$.
Note that $V^{\star}_1(\tilde \chi_i)=E(V^{\star s}(\tilde\chi_i,\tilde\chi_j,\tilde\chi_t)\mid \tilde\chi_i)=\frac{1}{3}u_iH_4(W_i)\neq 0$ with $H_4(W_i)=\frac{1}{2}E\left[(\eta_j\bar w_{ti}
+\eta_t\bar w_{ji})\mathcal{B}_{tj}\mid \tilde\chi_i\right]$. We can have %Hence $V_{n}^\star$ is not degenerate. Given Assumptions 4 and 5, under the standard U-statistic theory from \cite{serfling1980}, we can derive
\begin{eqnarray}
 V_n^{\star}=\frac{1}{n}\mathop{\sum}_{i=1 }^n u_iH_4(W_i)+o_p\left(\frac{1}
{\sqrt n}\right).
\end{eqnarray}
 similarly %as the discussion of $I_{n332}$, we have
 $I_{423}=O_p\left(\frac{1}
{\sqrt n^2}\right)=o_p\left(\frac{1}
{\sqrt n}\right)$.
Altogether, Proposition~4 is proved. \begin{flushright}
 $\Box$
 \end{flushright}

{\textbf{Proof of Proposition~5.}}

Recall that $\hat e_j-\tilde e_j=(V_j(\hat r)-V_j(r_0))Y_j-\sum_{t\neq j}w_{jt}(V_t(\hat r)-V_t(r_0))Y_t$ and $\overline Y_i-\sum_{k\neq i}w_{ik}\overline Y_t=\eta_i+G_i-\sum_{k\neq i}w_{ik}\bar Y_k$.
Again decompose $I_{n5}$ as:
\begin{eqnarray*}
&&I_{n5}=(\hat r-r_0)^{\top}\frac{1}{n(n-1)}\sum_{i=1}^n\sum_{j\neq i}  (\overline Y_i-\sum_{t\neq i}w_{it}\overline Y_t)(\overline Y_j-\sum_{k\neq j}w_{jk}\overline Y_k)^{\top}(\hat r-r_0)\mathcal{B}_{ij}+R_n\\
&&=(\hat r-r_0)^{\top}\frac{1}{n(n-1)}\sum_{i=1}^n\sum_{j\neq i}  (G_i-\sum_{t\neq i}w_{it}\overline Y_t)(G_j-\sum_{k\neq j}w_{jk}\overline Y_k)^{\top}\mathcal{B}_{ij}(\hat r-r_0)\\
&& \quad+(\hat r-r_0)^{\top}\frac{2}{n(n-1)}\sum_{i=1}^n\sum_{j\neq i}  \eta_i(G_j-\sum_{k\neq j}w_{jk}\overline Y_k)^{\top}\mathcal{B}_{ij}(\hat r-r_0)\\
&&\quad+(\hat r-r_0)^{\top}\frac{1}{n(n-1)}\sum_{i=1}^n\sum_{j\neq i}  \eta_i\eta_j^{\top}\mathcal{B}_{ij}(\hat r-r_0)+R_n\\
&&:=I_{n51}+I_{n52}+I_{n53}+R_n.
\end{eqnarray*}
Here $R_n$ is a higher order term.
By $\hat r-r_0=O_p\left(\frac{1}{\sqrt{n}}\right)$, $\mathop{\sup}_{X_i\in \Omega}|G_i-\sum_{k\neq i}w_{it}\overline Y_k|=O_p\left(h^s+\sqrt{\frac{\log n}{nh^p}}\right)$ and assumption 7, it can be easily obtained that $I_{n51}+I_{n52}=o_p\left(\frac{1}
{ n}\right)$, and the higher order term $R_n=o_p\left(\frac{1}
{ n}\right)$. Also note that $\left[\frac{1}{n(n-1)}
\sum_{i=1}^n\sum_{j\neq i}\eta_i\eta_j^{\top}
{\mathcal{B}}_{ij}\right]-
 E(\eta_1\eta_2^{\top}\mathcal{B}_{12})=o_p(1).$
Thus we can finish the proof of Proposition~5.
\begin{flushright}
$\Box$
\end{flushright}

\bigskip
\bigskip

{\textbf{Proof of Proposition~6.}}
Note that $e_i=a_nH(W_i)+u_i$ with $E(u_i\mid W)=0$. Thus we can rewrite $J_{n0}$ as
\begin{eqnarray}
J_{n0}&&=\frac{1}{n(n-1)}\sum_{i=1}^n\sum_{j\neq i}  e_i e_j\mathcal{B}_{ij}\nonumber\\
&&=\frac{1}{n(n-1)}\sum_{i=1}^n\sum_{j\neq i}  u_i u_j\mathcal{B}_{ij}+
2a_n\frac{1}{n(n-1)}\sum_{i=1}^n\sum_{j\neq i}  u_i H_j\mathcal{B}_{ij}+a_n^2\frac{1}{n(n-1)}\sum_{i=1}^n\sum_{j\neq i}  H_i H_j\mathcal{B}_{ij}\nonumber\\
&&:=J_{n01}+2a_nJ_{n02}+a_n^2J_{n03}
\end{eqnarray}
Obviously, $J_{n01}$ is a  degenerate U-statistic and thus $J_{n01}=O_p\left(\frac{1}{n}\right)$. Given Assumptions $4,5,9$,  U-statistic theory yields that
\begin{eqnarray}
&&J_{n02}=\frac{1}{n}\sum_{i=1}^nu_i H_{00}(W_i)+o_p\left(\frac{1}{\sqrt n}\right),\\
&&J_{n03}- \mu_0=o_p(1),
\end{eqnarray}
where $\mu_0=E(H_1H_2\mathcal{B}_{12})$ and $ H_{00}(W_1)=E(H_2\mathcal{B}_{12}\mid W_1).$
 Therefore $J_{n0}=O_p\left(\frac{1}{n}\right)+O_p\left(\frac{a_n}{\sqrt n}\right)+a_n^2\mu_0$. Thus, when $a_n$ is fixed, $J_{n0}=a_n^2\mu_0+o_p(1)$ and when $\sqrt{n}a_n\rightarrow \infty$, $J_{n0}/a_n^2=O_p\left(\frac{1}{na_n^2}\right)+O_p\left(\frac{1}{\sqrt na_n}\right)+ \mu_0=\mu_0+o_p(1).$ Further if $a_n=\frac{1}{\sqrt n}$, $nJ_{n0}=O_p(1)+O_p(1)+\mu_0$.
Proposition~6 is proved.
 \begin{flushright}
 $\Box$
 \end{flushright}

{\textbf{Proof of Proposition~7.}}
Under the alternative hypothesis, $e_i=a_nH(W_i)+u_i$ and then $J_{n1}$ can be decomposed as
\begin{eqnarray*}
&&J_{n1}=\frac{2}{n(n-1)}\sum_{i=1}^n\sum_{j\neq i} e_i(\tilde e_j-e_j) \mathcal{B}_{ij}\\
&&=\frac{2}{n(n-1)}\sum_{i=1}^n\sum_{j\neq i} u_i(\tilde e_j-e_j) \mathcal{B}_{ij}+\frac{2a_ n}{n(n-1)}\sum_{i=1}^n\sum_{j\neq i} H_i(\tilde e_j-e_j) \mathcal{B}_{ij}.
\end{eqnarray*}
Note that $\tilde e_j- e_j=(Y_j^*-\sum_{t\neq j}w_{jt} Y_t^*)-(Y_j^*-g(X_j))=g(X_j)-\sum_{t\neq j}w_{jt}g(X_t)-\sum_{t\neq j}w_{jt}u_t-a_n\sum_{t\neq j}w_{jt}H_t.$ Here $H_j=H(W_j)$. $J_{n1}$ can be further rewritten as
\begin{eqnarray*}
&&J_{n1}=\frac{2}{n(n-1)}\sum_{i=1}^n\sum_{j\neq i} u_i[g(X_j)-\sum_{t\neq j}w_{jt}g(X_t)-\sum_{t\neq j}w_{jt}u_t]\mathcal{B}_{ij} -a_n\frac{2}{n(n-1)^2}\sum_{i=1}^n\sum_{j\neq i}\sum_{t\neq j}u_i\tilde w_{jt} H_t \mathcal{B}_{ij}\nonumber\\
&&+a_n\frac{2}{n(n-1)}\sum_{i=1}^n\sum_{j\neq i} H_i[g(X_j)-\sum_{t\neq j}w_{jt}g(X_t)-\sum_{t\neq j}w_{jt}u_t]\mathcal{B}_{ij} -a_n^2\frac{2}{n(n-1)^2}\sum_{i=1}^n\sum_{j\neq i}\sum_{t\neq j}H_i\tilde w_{jt} H_t \mathcal{B}_{ij}\nonumber\\
&&:=2J_{n10}-2a_nJ_{n11}+2a_nJ_{n12}-2a_n^2J_{n13}.
\end{eqnarray*}
Obviously, $J_{n10}$ has a same asymptotic behavior as $I_{n1}$ and then $J_{n10}=-\frac{2}{n(n-1)}
\sum_{i=1}^n\sum_{j\neq i}u_iu_j\tilde{\mathcal{B}}_{ij}-\frac{1}{n}E[(u_1^2\bar w_{12}+u_2^2\bar w_{21})
\mathcal{B}_{12}]+o_p\left(\frac{1}{n}\right)$. Note that $E(\tilde w_{jt}H_t\mid X_t,W_j,W_i)=0$, and then similarly as the discussion on $I_{n33}$ in Proposition~3, we can derive that $a_nJ_{n11}=O_p\left(\frac{a_n}{n}\right)=o_p\left( a_n^2\right)$.

% thus if $a_n\rightarrow 0$, $ J_{n11}/a_n=o_p\left(\frac{1}{n}\right)$ and if $a_n>0$ is fixed, $J_{n11}=O_p\left(\frac{1}{n}\right)\rightarrow 0.$

Consider $J_{n12}$. Similarly as the proof for  $I_{n42}$ in Proposition~4, we can derive that
\begin{eqnarray}
J_{n12}=-\frac{1}{n}\sum_{i=1}^nu_iH_{10}(W_i)+o_p\left(\frac{1}{\sqrt n}\right),
\end{eqnarray}
where $\bar \chi_i=(W_i,u_i)$ and $H_{10}(W_i)=\frac{1}{2}E\left[(H_j\bar w_{ti}+H_t\bar w_{ji})\mathcal{B}_{tj}\mid \bar\chi_i\right]$. Thus we have $a_nJ_{n12}=O_p\left(\frac{a_n}{\sqrt n}\right)$.

Further, under the assumption $E(H_i\mid X_i)=0$ and assumption~3, we have $E(H_i\tilde w_{jt}H_t\mathcal{B}_{ij})=E[H_i\tilde w_{jt}\mathcal{B}_{ij}E(H_t\mid X_t,W_i,W_j)]=0$. Thus $J_{n13}=o_p(1)$ and $a_n^2J_{n13}=o_p(a_n^2).$

Thus, $J_{n1}=O_p\left(\frac{1}{ n}\right)+O_p\left(\frac{a_n}{\sqrt n}\right)+o_p(a_n^2)$. When $a_n$ is fixed, we have $J_{n1}\rightarrow 0$; when $\sqrt{n}a_n\rightarrow \infty$, $J_{n1}/a_n^2\rightarrow 0$; and when $a_n=\frac{1}{\sqrt n}$, $nJ_{n1}=nJ_{n10}-\frac{2}{\sqrt{n}}\sum_{i=1}^nu_iH_{10}(W_i)+o_p\left(1\right)$.
The proof of Proposition~7 is completed.
 \begin{flushright}
 $\Box$
 \end{flushright}

{\textbf{Proof of Proposition~8.}}
Recall that $\tilde e_j- e_j=(Y_j^*-\sum_{t\neq j}w_{jt} Y_t^*)-(Y_j^*-g(X_j))=g(X_j)-\sum_{t\neq j}w_{jt}g(X_t)-\sum_{t\neq j}w_{jt}u_t-a_n\sum_{t\neq j}w_{jt}H_t.$ Thus we have
\begin{eqnarray}
&&J_{n2}=\frac{1}{n(n-1)}\sum_{i=1}^n\sum_{j\neq i}  (\tilde e_i-e_i)(\tilde e_j- e_j)\mathcal{B}_{ij}\nonumber\\
&&=\frac{1}{n(n-1)}\sum_{i=1}^n\sum_{j\neq i}  [g(X_i)-\sum_{k\neq i}w_{ik}g(X_k)-\sum_{k\neq i}w_{ik}u_k][g(X_j)-\sum_{t\neq j}w_{jt}g(X_t)-\sum_{t\neq j}w_{jt}u_t]\mathcal{B}_{ij}+J_{n22}\nonumber\\
&&:=J_{n21}+J_{n22}
\end{eqnarray}
where
\begin{eqnarray}
&&J_{n22}= -2a_n\frac{1}{n(n-1)}\sum_{i=1}^n\sum_{j\neq i}  [g(X_i)-\sum_{k\neq i}w_{ik}g(X_k) ][\sum_{t\neq j}w_{jt}H_t]\mathcal{B}_{ij}\nonumber\\
&& \qquad \quad   +2a_n\frac{1}{n(n-1)}\sum_{i=1}^n\sum_{j\neq i}  [ \sum_{k\neq i}w_{ik}u_k][\sum_{t\neq j}w_{jt}H_t]\mathcal{B}_{ij} \nonumber\\ 
&&\qquad\quad +a_n^2\frac{1}{n(n-1)}\sum_{i=1}^n\sum_{j\neq i}  [\sum_{k\neq i}w_{ik}H_k][\sum_{t\neq j}w_{jt}H_t]\mathcal{B}_{ij}.\nonumber\\
&&\qquad :=-2a_nJ_{n221}+2a_nJ_{n222}+a_n^2J_{n223}.
\end{eqnarray}
$J_{n21}$ has the same asymptotic behavior as $I_{n2}$ in Proposition~2  to have  $J_{n21}=O_p\left(\frac{1}{ n}\right)$.

By the similar discussions as before, we can have that $J_{n221}=O_p(h^s)=o_p\left(\frac{1}{ \sqrt n}\right)$ and $J_{n223}=o_p(1)$. As for $J_{n222}$, noting that $E(u_k\mid W_1,\cdots,W_n)=0$ and $E(w_{jt}H_t\mid W_j,W_i)=0$ when $t\neq j,i$. Thus
the similar argument for dealing with $I_{n23}$ in Proposition~2 yields that  $J_{n222}=O_p\left(\frac{1}{ n}\right).$ Hence we  obtain that $J_{n2}=O_p\left(\frac{1}{ n}\right)-o_p\left(\frac{a_n}{ \sqrt n}\right)+O_p\left(\frac{a_n}{ n}\right)+o_p(a_n^2)$.
 When $a_n>0$ is fixed,  $J_{n2}=o_p(1)$; when $\sqrt{n}a_n\rightarrow\infty$, $J_{n2}/a_n^2=o_p(1)$; and when $a_n=\frac{1}{\sqrt n}$, $nJ_{n2}=O_p(1)=nJ_{n21}+o_p(1)$. Thus,  Proposition~8 is proved.
 \begin{flushright}
 $\Box$
 \end{flushright}

{\textbf{Proof of Proposition~9.}}
Decompose  $J_{n3}$  as
\begin{eqnarray}
&&J_{n3}=\frac{2}{n(n-1)}\sum_{i=1}^n\sum_{j\neq i} u_i(\hat e_j-\tilde e_j) \mathcal{B}_{ij}+\frac{2a_n}{n(n-1)}\sum_{i=1}^n\sum_{j\neq i} H_i(\hat e_j-\tilde e_j) \mathcal{B}_{ij}\nonumber\\
&&:=J_{n30}+J_{n31}.
\end{eqnarray}

Since $(\hat e_j-\tilde e_j)$ is not related to the alternative model, the asymptotic behavior of $J_{n30}$ is the same as $I_{n3}$ in Proposition~3, that is, $J_{n30}=O_p\left(\frac{1}{n}\right)=(\hat r-r_0)^{\top}\frac{2}{n}\sum_{i=1}^n u_iH_{3}(W_i)+o_p\left(\frac{1}{n}\right)$.

For $J_{n31}$, recall that
 $\hat e_j-\tilde e_j=(\hat r-r_0)^{\top}(\eta_j+G_j-\sum_{t\neq j} w_{jt}\overline Y_t)+o_p(\hat r-r_0)$. Then we have
\begin{eqnarray*}
&&J_{n31}=2a_n(\hat r-r_0)^{\top}\frac{1}{n(n-1)}\sum_{i=1}^n\sum_{j\neq i} H_i\eta_j \mathcal{B}_{ij}\\
&&\qquad +2a_n(\hat r-r_0)^{\top}\frac{1}{n(n-1)}\sum_{i=1}^n\sum_{j\neq i} H_i(G_j-\sum_{t\neq j} w_{jt}\overline Y_t) \mathcal{B}_{ij}+R_n\\
&&\qquad :=2a_n(\hat r-r_0)^{\top}J_{n311}+J_{n312}+R_n.
\end{eqnarray*}
Here $R_n$ is a higher order term. By $\hat r-r_0=O_p\left(\frac{1}{\sqrt n}\right)$ and $\mathop{\sup}_{X_j\in \Omega}|G_j-\sum_{t\neq j}w_{jt}\overline Y_t|=O_p\left(h^s+\sqrt{\frac{\log n}{nh^p}}\right)$,  $J_{n312}=o_p\left(\frac{a_n}{\sqrt n}\right)$. Also, it is easy to derive that $J_{n311}- E(H_1\eta_2\mathcal{B}_{12})=o_p(1)$. Thus the higher order term $R_n=o_p\left(\frac{a_n}{\sqrt n}\right).$
In summary, we have $J_{n3}=O_p\left(\frac{1}{ n}\right)+O_p\left(\frac{a_n}{\sqrt n}\right)+o_p\left(\frac{a_n}{\sqrt n}\right).$ Thus, when $a_n>0$ is fixed,  $J_{n3}=o_p(1)$; when $\sqrt{n}a_n\rightarrow \infty$, $J_{n3}/a_n^2=o_p(1)$; and when $a_n=\frac{1}{\sqrt n}$, $J_{n3}=(\hat r-r_0)^{\top}\frac{2}{n}\sum_{i=1}^n u_iH_{3}(W_i)+\frac{2}{\sqrt{n}}(\hat r-r_0)^{\top}\alpha+o_p\left(\frac{1}{ n}\right)$ {with $\alpha=E(H_1\eta_2\mathcal{B}_{12})$}.
The proof of Proposition~9 is completed.
 \begin{flushright}
 $\Box$
 \end{flushright}
{\textbf{Proof of Proposition~10.}}
Note that
\begin{eqnarray*}
&&J_{n4}=\frac{2}{n(n-1)}\sum_{i=1}^n\sum_{j\neq i}  (\hat  e_i-\tilde e_i)(\tilde e_j- e_j)\mathcal{B}_{ij}\\
&&=\frac{2}{n(n-1)}\sum_{i=1}^n\sum_{j\neq i}  (\hat e_i-\tilde e_i)[g(X_j)-\sum_{t\neq j}w_{jt}g(X_t)-\sum_{t\neq j}w_{jt}u_t ]\mathcal{B}_{ij}\nonumber\\
&&\quad -a_n\frac{2}{n(n-1)}\sum_{i=1}^n\sum_{j\neq i}\sum_{t\neq j}  (\hat e_i-\tilde e_i)w_{jt}H_t\mathcal{B}_{ij}\nonumber\\
&&:=J_{n41}-a_nJ_{n42}.
\end{eqnarray*}
The asymptotic behavior of $J_{n41}$ is the same as the one of $I_{n4}$ in Proposition~4 to have $J_{n41}=O_p\left(\frac{1}{ n}\right)=(\hat r-r_0)^{\top}\frac{1}{n}\sum_{i=1}^nu_iH_4(W_i)+o_p\left(\frac{1}{ n}\right)$.
As for $J_{n42}$, recall that
 $\hat e_j-\tilde e_j=(\hat r-r_0)^{\top}(\eta_j+G_j-\sum_{t\neq j} w_{jt}\overline Y_t)+o_p(\hat r-r_0)$, and then we can rewrite $J_{n42}$ as
 \begin{eqnarray*}
 &&J_{n42}=-a_n(\hat r-r_0)^{\top}\frac{2}{n(n-1)}\sum_{i=1}^n\sum_{j\neq i}\sum_{t\neq j}(G_j-\sum_{t\neq j} w_{jt}\overline Y_t)w_{jt}H_t\mathcal{B}_{ij}\\
 &&-a_n(\hat r-r_0)^{\top}\frac{2}{n(n-1)}\sum_{i=1}^n\sum_{j\neq i}\sum_{t\neq j}\eta_jw_{jt}H_t\mathcal{B}_{ij}+R_n\\
 &&:=J_{n421}-2a_n(\hat r-r_0)^{\top}J_{n422}+R_n.
 \end{eqnarray*}
Here $R_n$ is a higher order term. Very similarly as the proof for Proposition~9, %Noting that $\hat r-r_0=O_p\left(\frac{1}{\sqrt n}\right)$ and $\mathop{\sup}_{X_j\in \Omega}|G_j-\sum_{t\neq j}w_{jt}\overline Y_t|=O_p\left(h^s+\sqrt{\frac{\log n}{nh^p}}\right)$, hence it is easily obtain that $J_{n421}=o_p\left(\frac{a_n}{\sqrt n}\right)$.
%As for $J_{n422}$, we have
%\begin{eqnarray*}
%&&J_{n422}=\frac{1}{n(n-1)^2}\sum_{i=1}^n\sum_{j\neq i} \eta_j\tilde w_{ji}H_i\mathcal{B}_{ij}+\frac{1}{n(n-1)^2}\sum_{i=1}^n\sum_{j\neq i}\sum_{t\neq i, j}\eta_j\tilde w_{jt}H_t\mathcal{B}_{ij}\\
%&&:=J_{n4221}+J_{n4222}.
%\end{eqnarray*}
%Obviously, $J_{n4221}=O_p\left(\frac{1}{n}\right)$ and note that $E( \bar w_{jt}H_t\mid W_i,W_j)=0$, thus we have $J_{n422}=o_p(1)$. In summary,
we can also obtain that %$J_{n42}=o_p\left(\frac{a_n}{\sqrt n}\right)$ and
 $J_{n4}=(\hat r-r_0)^{\top}\frac{1}{n}\sum_{i=1}^nu_iH_4(W_i)+o_p\left(\frac{1}{ n}\right)+o_p\left(\frac{a_n}{\sqrt n}\right)$.
Therefore,  Proposition~10 is proved. %can be obtained.
  \begin{flushright}
 $\Box$
 \end{flushright}
{\textbf{Proof of Proposition~11.}}

Since $\hat e_i-\tilde e_i$ is not related to  the local alternative model, the asymptotic behavior of $J_{n5}$  is the same as that of $I_{n5}$ in Proposition~5 to have $J_{n5}=(\hat r-r_0)^{\top}A(\hat r -r_0)+o_p\left(\frac{1}{n}\right)$. Thus the proof of proposition~11 is completed.
  \begin{flushright}
 $\Box$
 \end{flushright}

\bibliography{test_ref}

\end{document}